\documentclass[12pt]{elsart}
\journal{Advances in Mathematics}
\usepackage{amscd,amssymb,verbatim}
\usepackage{amsmath}

\newcommand{\C}{{\mathbb C}}

\newcommand{\ch}{\operatorname{ch}}
\newcommand{\codim}{\operatorname{codim}}

\newcommand{\End}{\operatorname{End}}

\newcommand{\HH}{\operatorname{H}}

\newcommand{\Hol}{\operatorname{Hol}}
\newcommand{\Hom}{\operatorname{Hom}}

\newcommand{\Id}{\operatorname{Id}}

\newcommand{\Image}{\operatorname{Im}}

\newcommand{\Ind}{\operatorname{Ind}}

\newcommand{\KK}{\operatorname{K}}
\newcommand{\Ker}{\operatorname{Ker}}

\newcommand{\PHC}{\operatorname{PHC}}

\newcommand{\R}{{\mathbb R}}

\newcommand{\tr}{\operatorname{tr}}
\newcommand{\Tr}{\operatorname{Tr}}

\newcommand{\Z}{{\mathbb Z}}
\theoremstyle{plain}
\newtheorem{definition}{Definition}

\newtheorem{lemma}{Lemma}
\newtheorem{theorem}{Theorem}
\newtheorem{proposition}{Proposition}
\newtheorem{corollary}{Corollary}
\renewcommand{\rm}{\normalshape}
\begin{document}

\begin{frontmatter}

\title{Local Index Theory over Foliation Groupoids}

\author[gorokhovsky]{Alexander Gorokhovsky \thanksref{go}} and
\ead{gorokhov@euclid.colorado.edu}
\author[lott]{John Lott\corauthref{star} \thanksref{lo}}
\corauth[star]{Corresponding author.}
 \ead{lott@umich.edu}

\address[gorokhovsky]{Department of Mathematics, UCB 395,
University of Colorado, Boulder, CO~80309-0395, USA}

\address[lott]{Department of Mathematics,
University of Michigan, Ann Arbor, MI  48109-1043, USA}

\thanks[go]{The research
 was   supported in part by NSF grants DMS-0242780 and DMS-0400342.}

\thanks[lo]{The research was supported in part by
NSF grant DMS-0306242.}

\date{}
\begin{abstract}
We give a local proof of an index
theorem for a Dirac-type operator that is invariant
with respect to the action of a foliation groupoid $G$.
If $M$ denotes the space of units of $G$ then the input is a 
$G$-equivariant fiber bundle $P \rightarrow M$ along with
a $G$-invariant fiberwise Dirac-type operator $D$ on $P$. 
The index theorem is a formula for the pairing of the index
of $D$, as an element of a certain K-theory group, with a closed
graded trace on a certain noncommutative de Rham algebra
$\Omega^* {\mathcal B}$ associated to $G$. The proof is by means of
superconnections in the framework of noncommutative geometry.
\end{abstract}
\end{frontmatter}

\section{Introduction}

It has been clear for some time, especially since the work of
Connes \cite{Connes (1982)} and Renault \cite{Renault (1980)}, that
many interesting spaces in noncommutative geometry arise
from groupoids. For background information, we refer to
Connes' book \cite[Chapter II]{Connes (1994)}.
In particular, to a smooth groupoid $G$
one can assign its convolution algebra $C^\infty_c(G)$, which
represents a class of smooth functions on the noncommutative space
specified by $G$.

An important motivation for noncommutative geometry
comes from index theory.  The notion of groupoid allows one to
unify various index theorems that arise in the
literature, such as the Atiyah-Singer families index theorem
\cite{Atiyah-Singer (1971)}, the Connes-Skandalis foliation index
theorem \cite{Connes-Skandalis (1984)} and the Connes-Moscovici covering
space index theorem \cite{Connes-Moscovici (1990)}. All of these theorems
can be placed in the setting of
a proper cocompact action of a smooth groupoid $G$
on a manifold $P$. Given a $G$-invariant
Dirac-type operator $D$ on $P$, the construction of
\cite{Connes-Moscovici (1990)} allows one to form its analytic
index $\Ind_a$ as an element of the K-theory of the algebra
$C^\infty_c(G) \otimes {\mathcal R}$, where
${\mathcal R}$ is an algebra of infinite matrices whose entries decay
rapidly
\cite[Sections III.4, III.7.$\gamma$]{Connes (1994)}.
When composed with the trace on
${\mathcal R}$, the Chern character $\ch(\Ind_a)$
lies in the periodic cyclic homology group $\PHC_*(C^\infty_c(G))$. The index
theorem, at the level of Chern characters, equates
$\ch(\Ind_a)$ with a topological expression $\ch(\Ind_t)$.

We remark
that in the literature, one often sees the analytic index defined as an
element of K-theory of the groupoid $C^*$-algebra $C^*_r(G)$.
The index in $\KK_*(C^\infty_c(G) \otimes {\mathcal R})$ is a more
refined object.  However, to obtain geometric and topological
consequences from the index theorem, it appears that one has to
pass to $C^*_r(G)$; we refer to
\cite[Chapter III]{Connes (1994)} for discussion.  In this paper
we will work with $C^\infty_c(G)$.

We prove a local index theorem
for a Dirac-type operator that is  invariant
with respect to the action of a foliation groupoid.
In the terminology of Crainic-Moerdijk
\cite{Crainic-Moerdijk (2001)}, a foliation groupoid
is a smooth groupoid $G$ with discrete isotropy groups,
or equivalently, which is Morita equivalent to a
smooth \'etale groupoid.

A motivation for our work comes from the Connes-Skandalis index
theorem for a compact foliated manifold $(M, {\mathcal F})$ with a longitudinal
Dirac-type operator \cite{Connes-Skandalis (1984)}.
To a foliated manifold $(M, {\mathcal F})$ one can associate its
holonomy groupoid $G_{hol}$, which is an example of a foliation groupoid.
The general foliation index theorem equates $\Ind_a$ with a topological
index $\Ind_t$.  For details, we refer to
\cite[Sections I.5, II.8-9, III.6-7]{Connes (1994)}.

We now state the
index theorem that we prove.
Let $M$ be the space of units of a foliation groupoid $G$. It carries
a foliation ${\mathcal F}$.
Let $\rho$ be a closed holonomy-invariant
transverse current on $M$.
There is a corresponding universal class
$\omega_\rho \in \HH^*(BG; o)$, where $o$ is a certain orientation
character on the classifying space $BG$.
Suppose that $G$ acts freely, properly and cocompactly
on a manifold $P$. In particular, there is a submersion
$\pi : P \rightarrow M$.
There is an induced foliation $\pi^* {\mathcal F}$ of $P$ with
the same codimension as ${\mathcal F}$, satisfying
$T\pi^* {\mathcal F} \: = \: (d\pi)^{-1} T{\mathcal F}$.
Let $g^{TZ}$ be a smooth
$G$-invariant vertical Riemannian metric on $P$. Suppose that the
vertical tangent bundle $TZ$ is even-dimensional and
has a $G$-invariant spin structure.
Let $S^Z$ be the corresponding vertical spinor bundle.
Let $\widetilde{V}$ be an auxiliary $G$-invariant Hermitian vector bundle on
$P$ with a $G$-invariant Hermitian connection.
Put ${E} \: = \: S^Z \: \widehat{\otimes} \: \widetilde{V}$,
a $G$-invariant $\Z_2$-graded
Clifford bundle on $P$ which has a
$G$-invariant connection. The Dirac-type operator $Q$ acts
fiberwise on sections
of $E$. Let $D$ be its restriction to the sections of positive parity.
(The case of general $G$-invariant Clifford
bundles ${E}$ is completely analogous.) Let $\mu \: : \: P \rightarrow
P/G$ be the quotient map.
Then $P/G$
is a smooth compact manifold with a foliation $F \: = \: (\pi^*
{\mathcal F})/G$ satisfying
$(d\mu)^{-1} TF \: =  \: T \pi^*
{\mathcal F}$.  Put $V \: = \: \widetilde{V}/G$, a Hermitian
vector bundle on $P/G$ with a Hermitian connection $\nabla^V$.
The $G$-action on $P$ is classified by a map
$\nu : P/G \rightarrow BG$, defined up to homotopy.

The main point of this paper is to give a local proof of the
following theorem.
\begin{theorem} \label{first}
\begin{equation} \label{indexthm}
\langle \ch(\Ind D), \rho \rangle\: = \: \int_{P/G}
\widehat{A}({TF}) \: \ch(V) \: \nu^* \omega_\rho.
\end{equation}
\end{theorem}

Here $\Ind D$ lies in $\KK_*(C^\infty_c(G) \otimes {\mathcal R})$.
If $M$ is a compact foliated manifold and one takes
$P \: = \: G \: = \: G_{hol}$ then one recovers the
result of pairing the Connes-Skandalis theorem with $\rho$; see also
Nistor \cite{Nistor (1996)}.

In saying that we give a local proof of Theorem
\ref{first}, the word
``local'' is in the sense
of Bismut's proof of the Atiyah-Singer family index theorem
\cite{Bismut (1985)}.
In our previous paper
\cite{Gorokhovsky-Lott (2003)}
we gave a local proof of such a theorem in the \'etale case.
One can reduce Theorem \ref{first} to the \'etale case
by choosing a complete transversal $T$, i.e.
a submanifold of $M$, possibly disconnected, with $\dim(T) =
\codim({\mathcal F})$ and which intersects each
leaf of the foliation. Using $T$, one can reduce
the holonomy groupoid $G$ to a Morita-equivalent \'etale groupoid
$G_{et}$. We gave a
local proof of
Connes' index theorem concerning an \'etale groupoid $G_{et}$ acting
freely,
properly and cocompactly on a manifold $P$, preserving a fiberwise
Dirac-type operator $Q$ on $P$.
Our local proof has since been used by Leichtnam and Piazza
to prove an index theorem for foliated manifolds-with-boundary
\cite{Leichtnam-Piazza (2003)}.

In the present paper we give a local proof of
Theorem \ref{first} working directly with foliation groupoids.
In particular, the
new proof avoids the noncanonical choice of a complete transversal $T$.

The overall method of proof is by means of superconnections in the
context of noncommutative geometry, as in
\cite{Gorokhovsky-Lott (2003)}. However, there are conceptual
differences with respect to \cite{Gorokhovsky-Lott (2003)}.
As in \cite{Gorokhovsky-Lott (2003)}, we first establish an
appropriate differential calculus on the noncommutative space determined
by a
foliation groupoid $G$. The notion of ``smooth functions'' on the
noncommutative space is clear, and is given by the elements of the
convolution algebra
${\mathcal B} = C^\infty_c(G)$.
We define a certain graded
algebra $\Omega^*{\mathcal B}$ which plays the role of the differential
forms on the noncommutative space.
The algebra $\Omega^*\mathcal {B}$ is equipped with a degree-$1$
derivation $d$, which is the analog of the de Rham differential.
Unlike in the \'etale case, it turns out that in general, $d^2 \neq 0$.
The reason for this is that to define
$d$, we must choose a horizontal distribution $T^H M$ on $M$, where
``horizontal'' means transverse to
${\mathcal F}$. In general $T^H M$ is not integrable,
which leads to the nonvanishing of $d^2$. This issue does not arise
in the \'etale case.

As we wish to deal with superconnections in such a context, we must first
understand how to do Chern-Weil theory when $d^2 \neq 0$. If
$d^2$ is given by commutation with a $2$-form then a trick of Connes
\cite[Chapter III.3, Lemma 9]{Connes (1994)}
allows one to construct a new complex
with $d^2 = 0$, thereby reducing to the usual case.
We give a somewhat more general formalism that may be useful in other
contexts.  It assumes that
for the relevant
${\mathcal B}$-module $\mathcal{E}$ and connection
$\nabla \: : \: \mathcal{E} \rightarrow \Omega^1 \mathcal{B} 
\otimes_{\mathcal{B}}
\mathcal{E}$, there is a linear map
$l : \mathcal{E} \to \Omega^2 \mathcal{B} \otimes_{\mathcal{B}}
\mathcal{E}$ such that
\begin{equation} \label{l}
l(b \xi) -\: b \: l(\xi) \: = \: d^2(b) \: \xi
\end{equation}
and
\begin{equation} \label{ll}
l(\nabla \xi) \: = \: \nabla l(\xi)
\end{equation}
for $b \in \mathcal{B}$, $\xi \in \mathcal{E}$. With this additional
structure, we show in Section \ref{chern} how to do Chern-Weil theory,
both for connections and superconnections on a ${\mathcal B}$-module
$\mathcal{E}$. In the case when $d^2$ is a commutator, one recovers
Connes' construction of Chern classes.

Next, we consider certain ``homology classes'' of the noncommutative space.
A graded trace on $\Omega^*{\mathcal B}$ is said to be
closed if it annihilates $\Image(d)$. A closed holonomy-invariant transverse
current $\rho$ on the space of units $M$ gives a closed graded trace
on $\Omega^* {\mathcal B}$.

The action of $G$ on $P$ gives rise to
a left ${\mathcal B}$-module
${\mathcal E}$,
which essentially consists of compactly-supported sections of $E$
coupled to a vertical density. We extend $\mathcal{E}$ to a
left-$\Omega^*\mathcal{B}$ module $\Omega^*\mathcal{E}$
of ``$\mathcal{E}$-valued differential forms''.
There is a natural linear map $l : \mathcal{E} \to \Omega^2\mathcal{E}$
satisfying (\ref{l}) and (\ref{ll}).

We then consider
the Bismut superconnection $A_s$
on $\mathcal{E}$. The formal expression for its Chern character
involves $e^{- \: A_s^2 \: + \: l}$. The latter is well-defined in
$\Hom^\omega ({\mathcal E}, \Omega^* {\mathcal E})$, an algebra
consisting of
rapid-decay kernels.
We construct a graded trace $\tau :
\Hom^\omega ({\mathcal E}, \Omega^* {\mathcal E}) \rightarrow
\Omega^* \mathcal{B}$.  This allows us to define the Chern
character of the superconnection by
\begin{equation}
\ch(A_s) \: = \: {\mathcal R} \left(\tau e^{- \: A_s^2 \: + \: l}
\right).
\end{equation}
Here ${\mathcal R}$ is the rescaling operator which, for $p$ even,
multiplies a $p$-form by $(2 \pi i)^{- \: \frac{p}{2}}$.

Now let $\rho$ be a closed holonomy-invariant transverse current on $M$
as above. Then
$\rho(\ch(A_s))$ is defined and we compute its  limit
when $s \to 0$, to obtain a differential form version
of the right-hand-side of (\ref{indexthm}).  (In the case when
$P \: = \: G \: = \: G_{hol}$ an
analogous computation was done by Heitsch
\cite[Theorem 2.1]{Heitsch (1995)}).

Next, we use the argument of
\cite[Section 5]{Gorokhovsky-Lott (2003)} to show that
for all $s > 0$,
$\langle \ch(\Ind D), \rho \rangle\ = \rho(\ch(A_s))$.
(In the case when $P \: = \: G \: = \: G_{hol}$, this was
shown under some further restrictions by Heitsch
\cite[Theorem 4.6]{Heitsch (1995)} and Heitsch-Lazarov
\cite[Theorem 5]{Heitsch-Lazarov (1999)}.)
This proves Theorem \ref{first}.

We note that our extension of \cite{Gorokhovsky-Lott (2003)}
from \'etale groupoids to foliation groupoids
is only partial. The local index theorem of
\cite{Gorokhovsky-Lott (2003)} allows for pairing with more general
objects than transverse currents, such as
the Godbillon-Vey class.  The paper \cite{Gorokhovsky-Lott (2003)} used
a bicomplex $\Omega^{*,*} {\mathcal B}$ of forms, in which the
second component consists of forms in the ``noncommutative'' direction.
There was also a connection $\nabla$ on ${\mathcal E}$ which involved
a differentiation in the noncommutative direction.  In the setting of
a foliation groupoid, one again has a bicomplex $\Omega^{*,*} {\mathcal B}$
and a connection $\nabla$. However, (\ref{ll}) is not satisfied.
Because of this we work instead with the smaller complex of forms
$\Omega^{*,0} {\mathcal B}$, where this problem does not arise.

The paper is organized as follows. In Section \ref{chern} we
discuss Chern-Weil theory in the context of a graded algebra
with derivation whose square is nonzero. In Section
\ref{differential} we describe the differential algebra
$\Omega^* {\mathcal B}$
associated to a foliation
groupoid $G$. In Section \ref{superconnection} we add a manifold $P$ on which
$G$ acts properly.
We define a certain left-${\mathcal B}$
module ${\mathcal E}$ and superconnection
$A_s$ on ${\mathcal E}$. We compute the $s \rightarrow 0$ limit of
$\rho(\ch(A_s))$.
In Section \ref{Index Theorem} we explain the relation between the
superconnection computations and the K-theoretic index,
construct the cohomology class $\omega_{\rho} \in
\HH^*(BG; o)$ and prove Theorem \ref{first}. We show that Theorem 
\ref{first} implies some well-known index theorems.

In an appendix to this paper we give a technical improvement to our previous
paper \cite{Gorokhovsky-Lott (2003)}. The index theorem in
\cite{Gorokhovsky-Lott (2003)} assumed that the closed graded trace
$\eta$ on $\Omega^*(B, \C \Gamma)$ extended to an algebra of rapidly
decaying forms
$\Omega^*(B, {\mathcal B}^\omega)$. The appearance of
$\Omega^*(B, {\mathcal B}^\omega)$ was due to the noncompact support of the
heat kernel, which affects the trace of the superconnection Chern character.
In the appendix we show how to replace $\Omega^*(B, {\mathcal B}^\omega)$ by
$\Omega^*(B, \C \Gamma)$, by
using finite propagation speed methods.
Let $f \in
C^\infty_c(\R)$ be a smooth even function with support in
$[-\epsilon, \epsilon]$. Let $\widehat{f}$ be its Fourier transform.
We can define $\widehat{f}(A_s)$ and show that
$\eta \left( {\mathcal R} \: \tau \widehat{f}(A_s)
\right)$ is defined for graded traces
$\eta$ on $\Omega^*(B, \C \Gamma)$.
We prove the corresponding
analog of \cite[Theorem 3]{Gorokhovsky-Lott (2003)}, with the Gaussian
function in the definition of the Chern character
replaced by an appropriate function $\widehat{f}$. This then implies
the result stated in \cite[Theorem 3]{Gorokhovsky-Lott (2003)}
without the condition of $\eta$ being extendible to
$\Omega^*(B, {\mathcal B}^\omega)$. We remark that this issue of
replacing $\Omega^*(B, {\mathcal B}^\omega)$ by
$\Omega^*(B, \C \Gamma)$ does not arise in the present paper.

More detailed summaries are given at the beginnings of the sections.

We thank the referee for useful comments.

\section{The Chern Character}\label{chern}

In this section we collect some algebraic facts needed to define the
Chern character of a superconnection in our setting. We consider
an algebra ${\mathcal B}$ and
a graded algebra $\Omega^*$ with $\Omega^0 \: = \: {\mathcal B}$.
We assume that $\Omega^*$ is equipped
with a degree-$1$ derivation $d$ whose square may be nonzero.
If ${\mathcal E}$ is a left ${\mathcal B}$-module then the notion
of a connection $\nabla$ on ${\mathcal E}$ is the usual one from
noncommutative geometry; see Connes
\cite[Section III.3, Definition 5]{Connes (1994)} and Karoubi
\cite[Chapitre 1]{Karoubi (1987)}.
We assume the additional structure of a
map $l$ satisfying (\ref{l}) and (\ref{ll}).
We show that $\nabla^2 - l$ is then the
right notion of curvature. If ${\mathcal E}$ is a finitely-generated
projective ${\mathcal B}$-module then we carry out Chern-Weil theory
for the connection $\nabla$, and show how it extends to the
case of a superconnection $A$. Many of the lemmas in this section are
standard in the case when $d^2 = 0$ and $l = 0$, but we present them
in detail in order to make clear what goes through to the case when
$d^2 \neq 0$. In the case when $d^2$ is given by a commutator, the
Chern character turns out to be the same as what one would get using
Connes' $X$-trick
\cite[Section III.3, Lemma 9]{Connes (1994)}.

 Let ${\mathcal
B}$ be an algebra over $\C$, possibly nonunital. Let $\Omega \: =
\: \bigoplus_{i=1}^\infty \Omega^i$ be a graded algebra with
$\Omega^0 \: = \: {\mathcal B}$. Let $d \: : \: \Omega^*
\rightarrow \Omega^{*+1}$ be a graded derivation of $\Omega^*$.
Define $\alpha  \: : \: \Omega^* \rightarrow \Omega^{*+2}$ by
$\alpha \: = \: d^2$; then for all $\omega, \omega^\prime \in
\Omega^*$,
\begin{equation}
\alpha(d\omega) \: = \: d \alpha(\omega), \: \: \: \: \: \:
\alpha(\omega \omega^\prime) \: = \: \alpha(\omega) \: \omega^\prime
\: + \: \omega \: \alpha(\omega^\prime).
\end{equation}

By a graded trace, we will mean a linear functional
$\eta \: : \: \Omega^* \rightarrow \C$ such that
\begin{equation}
\eta(\alpha(\omega)) \: = \: 0, \: \:  \: \:  \: \:  \: \:  \: \:
\: \: \eta([\omega, \omega^\prime]) \: = \: 0
\end{equation}
for all $\omega, \omega^\prime \in \Omega^*$. Define
$d^t \eta$ by $(d^t \eta)(\omega) \: = \: \eta(d\omega)$. Then the
graded traces on $\Omega^*$ form a complex with differential $d^t$.
A graded trace $\eta$ will be said to be closed if $d^t \eta \: = \: 0$,
i.e. for all $\omega \in \Omega^*$, $\eta(d \omega) \: = \: 0$. \\ \\
{\bf Example 1 : } Let $E$ be a complex vector bundle over a smooth
manifold $M$. Let $\nabla^E$ be a connection on $E$, with curvature
$\theta^E \in
\Omega^2(M; \End(E))$. Put ${\mathcal B} \: = \: C^\infty(M; \End(E))$
and $\Omega^* \: = \: \Omega^*(M; \End(E))$. Let $d$ be the extension of
the connection $\nabla^E$ to $\Omega^*(M; \End(E))$. Then $\alpha(\omega)
\: =
\: \theta^E \: \omega \: - \: \omega \: \theta^E$. If
$c$ is a closed current on $M$ then we obtain a closed graded trace $\eta$ on
$\Omega^*$ by $\eta(\omega) \: = \: \int_c \tr(\omega)$. \\

Let ${\mathcal E}$ be a left ${\mathcal B}$-module.
We assume that there is a
$\C$-linear
map $l \: : \: {\mathcal E} \rightarrow \Omega^2 \otimes_{\mathcal B} {\mathcal E}$
such that for all $b \in {\mathcal B}$ and $\xi \in {\mathcal E}$,
\begin{equation}
l(b \xi) \: = \: \alpha(b) \: \xi \: + \: b \: l(\xi).
\end{equation}
{\bf Example 2 : } Suppose that for some $\theta \in \Omega^2$,
$\alpha( \omega) \: = \: \theta \:
\omega \: - \: \omega \theta$. Then
we can take $l(\xi) \: = \: \theta \xi$.
\begin{lemma}\label{lemma1}
There is an extension of $l$ to a linear
map $l \: : \: \Omega^* \otimes_{\mathcal B} {\mathcal E}
\rightarrow \Omega^{*+2} \otimes_{\mathcal B} {\mathcal E}$ so that
for $\omega \in \Omega^*$ and $\mu \in \Omega^* \otimes_{\mathcal B} {\mathcal
E}$,
\begin{equation}
l(\omega \mu) \: = \: \alpha(\omega) \: \mu \: + \: \omega \: l(\mu).
\end{equation}
\end{lemma}
\begin
{pf}
We define $l \: : \: \Omega^* \otimes_{\C} {\mathcal E} \rightarrow
\Omega^{*+2} \otimes_{\mathcal B} {\mathcal E}$ by
\begin{equation}
l(\omega \otimes \xi) \: = \: \alpha(\omega) \: \xi \: + \: \omega \:
l(\xi).
\end{equation}
Then for $b \in {\mathcal B}$,
\begin{align}
l(\omega b \otimes \xi) \: & = \: \alpha(\omega b) \: \xi \: + \: \omega
b  \: l(\xi) \: = \: \alpha(\omega) \: b \: \xi \: + \: \omega \:
\alpha(b)
\: \xi
\: + \: \omega \: b \: l(\xi) \\
& = \: \alpha(\omega) \: b \xi \: + \: \omega \:
l(b \xi) \: = \: l(\omega \otimes b\xi). \notag
\end{align}
Thus $l$ is defined on $\Omega^* \otimes_{\mathcal B} {\mathcal E}$.
Next, for $\omega, \omega^\prime \in \Omega^*$ and $\xi \in {\mathcal E}$,
\begin{align}
l(\omega ( \omega^\prime \xi)) \: & = \:
\alpha(\omega \omega^\prime) \: \xi \: + \: \omega \omega^\prime \:
l(\xi) \: = \: \alpha(\omega) \:  \omega^\prime \: \xi \: + \:
\omega \: \alpha(\omega^\prime) \: \xi \: + \: \omega
\omega^\prime \: l(\xi) \\
& = \:
\alpha(\omega) \:  \omega^\prime \: \xi \: + \:
\omega \: l(\omega^\prime \: \xi). \notag
\end{align}
This proves the lemma.
\end{pf}

Let $\nabla \: : \: {\mathcal E} \rightarrow \Omega^1 \otimes_{\mathcal B}
{\mathcal E}$
be a connection, i.e. a $\C$-linear map satisfying
\begin{equation}
\nabla(b \xi) \: = \: db \otimes \xi \: + \: b \nabla \xi
\end{equation}
for all $b \in {\mathcal B}$, $\xi \in {\mathcal E}$. Extend $\nabla$ to a
$\C$-linear map
$\nabla \: : \: \Omega^* \otimes_{\mathcal B} {\mathcal E} \rightarrow
\Omega^{*+1} \otimes_{\mathcal B} {\mathcal E}$ so that for all
$\omega \in \Omega^*$ and $\xi \in {\mathcal E}$,
\begin{equation}
\nabla(\omega \xi) \: = \: d\omega \otimes \xi \: + \:
(-1)^{|\omega|} \: \omega \nabla \xi.
\end{equation}
We assume that for all
$\xi \in {\mathcal E}$,
\begin{equation} \label{comm}
l(\nabla \xi) \: = \: \nabla l(\xi).
\end{equation}

\begin{lemma}\label{lemma2}
$\nabla^2 \: - \: l \: : \:  {\mathcal E} \rightarrow
\Omega^{2} \otimes_{\mathcal B} {\mathcal E}$ is left-${\mathcal B}$-linear.
\end{lemma}
\begin{pf}
For $b \in {\mathcal B}$ and $\xi \in {\mathcal E}$,
\begin{align}
(\nabla^2 \: - \: l)(b \xi) \: & = \: \nabla(db \otimes \xi \: + \: b
\nabla \xi) \: - \: l(b \xi) \: = \: d^2 b \otimes \xi \:
+ \: \: b \nabla^2 \xi \: - \:
l(b \xi) \\
& = \:  \alpha(b) \xi \: + \: b \nabla^2 \xi \: - \:
l(b \xi) \: = \: b \: (\nabla^2 \: - \: l)(\xi). \notag
\end{align}
This proves the lemma.
\end{pf}

Put $\Omega^*_{ab} \: = \: \Omega^*/[\Omega^*, \Omega^*]$, the quotient
by the graded commutator, with the induced $d$.
For simplicity, in the rest of this section we assume that
${\mathcal B}$ is unital and
${\mathcal E}$ is a finitely-generated projective left ${\mathcal B}$-module.
Consider the graded algebra $\Hom_{\mathcal B} \left( {\mathcal E},
\Omega^* \otimes_{\mathcal B} {\mathcal E} \right) \: \cong \:
\End_{\Omega^*} \left( \Omega^* \otimes_{\mathcal B} {\mathcal E} \right)$.
There is a graded trace on $\Hom_{\mathcal B} \left( {\mathcal E},
\Omega^* \otimes_{\mathcal B} {\mathcal E} \right)$, with value in
$\Omega^*_{ab}$, defined as follows.
Write ${\mathcal E}$ as ${\mathcal B}^N e$ for some idempotent
$e \in M_N({\mathcal B})$. Then any
$T \in \Hom_{\mathcal B} \left( {\mathcal E},
\Omega^* \otimes_{\mathcal B} {\mathcal E} \right)$ can be represented as
right-multiplication on ${\mathcal B}^N e$ by a matrix
$T \in M_N( \Omega^*)$ satisfying $T \: = \: eT \: = \: Te$. By
definition $\tr(T) \: = \: \sum_{i=1}^N T_{ii} \mod{[\Omega^*, \Omega^*]}$.
It is independent of the representation of ${\mathcal E}$ as ${\mathcal B}^N e$.

Given $T_1, T_2 \in
\End_{\Omega^*} \left(
\Omega^* \otimes_{\mathcal B} {\mathcal E} \right)$, define their (graded)
commutator by
\begin{equation}\label{2.12}
[T_1, T_2] \: = \: T_1 \circ T_2 \: - \: (-1)^{|T_1||T_2|} \: T_2 \circ
T_1.
\end{equation}
For $T \in \End_{\Omega^*} \left(
\Omega^* \otimes_{\mathcal B} {\mathcal E} \right)$, define
$[\nabla, T] \in \End_{\C} \left(
\Omega^* \otimes_{\mathcal B} {\mathcal E} \right)$ by
\begin{equation}
[\nabla, T](\mu) \: = \: (-1)^{|\mu|} \left( \nabla(T (\mu)) \: - \:
T(\nabla \mu) \right)
\end{equation}
for $\mu \in \Omega^* \otimes_{\mathcal B} {\mathcal E}$.
\begin{lemma}
$[\nabla, T] \in \End_{\Omega^*} \left(
\Omega^* \otimes_{\mathcal B} {\mathcal E} \right)$.
\end{lemma}
\begin{pf}
Given $\omega \in \Omega^*$ and $\mu \in \Omega^* \otimes_{\mathcal B} {\mathcal E}$,
\begin{align}
[\nabla, T](\omega \mu) \: & = \: (-1)^{|\omega|+|\mu|} \left(
\nabla(T(\omega \mu)) \: - \: T(\nabla (\omega \mu)) \right) \\
& = \:
(-1)^{|\omega|+|\mu|} \left(
\nabla(\omega T(\mu)) \: - \: T((d\omega) \mu \: + \:
(-1)^\omega \omega \nabla \mu)
\right) \notag \\
& = \:
(-1)^{|\omega|+|\mu|} \left(
(d\omega) T(\mu) \: + \: (-1)^{\omega} \omega \nabla(T(\mu))
\: - \: (d\omega) T(\mu) \: - \: (-1)^{\omega} \: \omega T(\nabla \mu)
\right) \notag \\
& = \: \omega \: [\nabla, T](\mu). \notag
\end{align}
This proves the lemma.
\end{pf}
\begin{lemma}
Given $T_1, T_2 \in
\End_{\Omega^*} \left(
\Omega^* \otimes_{\mathcal B} {\mathcal E} \right)$,
\begin{equation}\label{2.15}
[\nabla, T_1 \circ T_2] \: = \: T_1 \circ [\nabla, T_2] \: + \:
(-1)^{|T_2|} \: [\nabla, T_1] \circ T_2.
\end{equation}
\end{lemma}
\begin{pf}
Given $\mu \in \Omega^* \otimes_{\mathcal B} {\mathcal E}$,
\begin{equation}
[\nabla, T_1 \circ T_2](\mu) \: = \: (-1)^{|\mu|}
\left\{ \nabla(T_1(T_2(\mu))) \: - \: T_1(T_2(\nabla (\mu))) \right\},
\end{equation}
\begin{equation}
\left( T_1 \circ [\nabla, T_2] \right)(\mu) \: = \: (-1)^{|\mu|} \:
T_1 \left( \nabla(T_2(\mu)) \: - \: T_2(\nabla (\mu)) \right)
\end{equation}
and
\begin{equation}
\left( [\nabla, T_1] \circ T_2 \right)(\mu) \: = \:
[\nabla, T_1](T_2(\mu)) \: = \:
(-1)^{|T_2(\mu)|}
\left\{ \nabla(T_1(T_2(\mu))) \: - \: T_1(\nabla (T_2(\mu))) \right\}.
\end{equation}
The lemma follows.
\end{pf}
\begin{lemma}
Given $T_1, T_2 \in
\End_{\Omega^*} \left(
\Omega^* \otimes_{\mathcal B} {\mathcal E} \right)$,
\begin{equation}
[\nabla, [T_1, T_2]] \: = \: [T_1, [\nabla, T_2]]
\: +
\: (-1)^{|T_2|} \: [[\nabla, T_1], T_2].
\end{equation}
\end{lemma}
\begin{pf}
This follows from (\ref{2.12}) and (\ref{2.15}).  We omit the
details.
\end{pf}
\begin{lemma}
For $T \in
\End_{\Omega^*} \left(
\Omega^* \otimes_{\mathcal B} {\mathcal E} \right)$,
\begin{equation}
\tr([\nabla, T]) \: = \: d \tr(T) \in \Omega^*_{ab}.
\end{equation}
\end{lemma}
\begin{pf}
Let us write ${\mathcal E} \: = \: {\mathcal B}^N e$ for an idempotent $e \in
M_N({\mathcal B})$.
Given $A \in \Hom_{\mathcal B} ({\mathcal E},
\Omega^1 \otimes_{\mathcal B} {\mathcal E})$, it acts on ${\mathcal B}^N e$ on
the right by a matrix
$A \in M_N(\Omega^1)$ with $A \: = \: eA \: = \: Ae$.
Then there is some $A \in \Hom_{\mathcal B} ({\mathcal E},
\Omega^1 \otimes_{\mathcal B} {\mathcal E})$ so that
for $\mu \in \Omega^* \otimes_{\mathcal B} {\mathcal E}
\: = \: (\Omega^*)^N e$,
\begin{equation}
\nabla(\mu) \: = \: (d\mu) \: e
\: + \: (-1)^{|\mu|} \: \mu \: A;
\end{equation}
in fact, this equation defines $A$.

An element $T \in \End_{\Omega^*} \left(
\Omega^* \otimes_{\mathcal B} {\mathcal E} \right)$ 
acts by right multiplication
on $\Omega^* \otimes_{\mathcal B} {\mathcal E}
\: = \: (\Omega^*)^N e$
by a matrix $T \in M_N(\Omega^*)$ satisfying $T \: = \: eT \: = \: Te$.
Then for $\xi \in {\mathcal E} \: = \: {\mathcal B}^N e$,
\begin{align}
[\nabla, T](\xi) \: & = \: \nabla(\xi T) \: - \: (\nabla(\xi)) T \: = \:
\left\{ d(\xi T) \: e \: + \: (-1)^{|T|} \: \xi T A \right\} \: - \:
\left\{ (d\xi) \: e \:  + \: \xi A \right\} T \\
& = \: \xi \left( (dT) \: e \: + \: (-1)^{|T|} \: T A \: - \:
 A T \right) \notag
\end{align}
Thus $[\nabla, T]$ acts as right multiplication by the matrix
\begin{equation}
e (dT) e \: + \: (-1)^{|T|} \: T A \: - \:
 A T,
\end{equation}
and so $\tr([\nabla, T]) \: \equiv \: \tr(e(dT)e)$. On the other hand,
using the identity $e(de) e \: = \: 0$ and taking the trace of $N \times
N$ matrices, we obtain
\begin{align}
d \tr(T) \: & = \: d \tr(eTe) \: = \:
\tr \left( (de) T e \: + \: e(dT)e \: + \: (-1)^{|T|} \: eT(de) \right)
\\
& = \:
\tr \left( (de)e T e \: + \: e(dT)e \: + \: (-1)^{|T|} \: eTe(de) \right)
\notag \\
& \equiv \:
\tr \left( e(de)e T \: + \: e(dT)e \: + \: (-1)^{|T|} \: Te(de)e \right)
\: = \:
\tr(e(dT) e). \notag
\end{align}
This proves the lemma.
\end{pf}
\begin{lemma}
$[\nabla, \nabla^2 - l] \: = \: 0$.
\end{lemma}
\begin{pf}
This follows from (\ref{comm}).
\end{pf}

\begin{definition}
The Chern character form of $\nabla$ is
\begin{equation}
\ch(\nabla) \: = \: \tr \left( e^{- \: \frac{\nabla^2 - l}{2\pi i}} \right) \in
\Omega^*_{ab}.
\end{equation}
\end{definition}

\begin{lemma}\label{lemma8}
Given ${\mathcal E}$, if
$\eta$ is a closed graded trace on $\Omega^*$
then $\eta(\ch(\nabla))$ is independent
of the choice of $\nabla$. If $\eta_1$ and $\eta_2$ are homologous closed
graded traces then $\eta_1(\ch(\nabla)) \: = \: \eta_2(\ch(\nabla))$.
\end{lemma}
\begin{pf}
Let $\nabla_1$ and $\nabla_2$ be two connections on ${\mathcal E}$.
For $t \in [0,1]$, define a
connection by $\nabla(t) \: = \: t \nabla_2 \: + \: (1-t) \nabla_1$.
Then $\frac{d\nabla}{dt} \: = \: \nabla_2 \: - \: \nabla_1 \in
\Hom_{\mathcal B} ({\mathcal E}, \Omega^1 \otimes_{\mathcal B} {\mathcal E})$.
We claim that $\eta(\ch(\nabla(t)))$ is independent of $t$.
As $\frac{d(\nabla^2 \: - \: l)}{dt}  \: = \: \nabla \frac{d\nabla}{dt}
\: + \: \frac{d\nabla}{dt} \nabla$, we have
\begin{align}
\frac{d\ch(\nabla)}{dt} \: & = \: - \: \frac{1}{2\pi i} \: \tr \left(
\left( \nabla \frac{d\nabla}{dt}
\: + \: \frac{d\nabla}{dt} \nabla \right) \:
e^{- \: \frac{\nabla^2 - l}{2\pi i}}
\right) \: = \:
- \: \frac{1}{2\pi i} \: \tr \left(
\left[\nabla, \frac{d\nabla}{dt}
\: e^{- \: \frac{\nabla^2 - l}{2\pi i}} \right]
\right) \\
& = \: - \: \frac{1}{2\pi i} \: d \: \tr \left(
\frac{d\nabla}{dt} \: e^{- \: \frac{\nabla^2 - l}{2\pi i}} \right). \notag
\end{align}
Then
\begin{equation}\label{2.28}
\ch(\nabla_2) \: - \: \ch(\nabla_1) \: = \: - \:
\frac{1}{2\pi i} \: d \: \int_0^1 \tr \left(
(\nabla_2 \: - \: \nabla_1) \:
e^{- \: \frac{\nabla(t)^2 - l}{2\pi i}} \right) \: dt,
\end{equation}
from which the claim follows. We note after expanding the
exponential in (\ref{2.28}), the integral gives an expression that
is purely algebraic in $\nabla_1$ and $\nabla_2$.

If $\eta_1$ and $\eta_2$ are homologous then there is a graded
trace $\eta^\prime$ such that $\eta_1 \: - \: \eta_2 \: = \: d^t \eta^\prime$.
Thus
\begin{equation}
\eta_1(\ch(\nabla)) \: - \: \eta_2(\ch(\nabla)) \: = \:
\eta^\prime (d \ch(\nabla)).
\end{equation}
However,
\begin{equation}
d\ch(\nabla) \: = \: d \tr \left( e^{- \: \frac{\nabla^2 - l}{2\pi i}} \right)
\: = \: \tr \left( \left[ \nabla, e^{- \: \frac{\nabla^2 - l}{2\pi i}} \right]
\right) \: = \: 0.
\end{equation}
This proves the lemma.
\end{pf}
\noindent
{\bf Example 3 : } With the notation of Example 1, let $F$ be another
complex vector bundle on $M$, with connection $\nabla^F$. Put ${\mathcal E}
\: = \: C^\infty(M; E \otimes F)$, with $l(\xi) \: = \: (\theta^E \otimes
I) \: \xi$ for
$\xi \in {\mathcal E}$. Let
$\nabla$ be the tensor product of $\nabla^E$ and $\nabla^F$.
Then one finds that $\eta(\ch(\nabla)) \: = \:
\int_c \ch(\nabla^F)$. \\

If ${\mathcal E}$ is $\Z_2$-graded,
let
$A \: : \:
{\mathcal E} \rightarrow \Omega^* \otimes_{\mathcal B} {\mathcal E}$ be a
superconnection. Then there are obvious extensions of the results of this
section.  In particular, let ${\mathcal R}$ be the rescaling operator on
$\Omega^{even}_{ab}$ which multiplies an element of $\Omega^{2k}_{ab}$ by
$(2 \pi i)^{-k}$.

\begin{definition}
The Chern character form of $A$ is
\begin{equation}
\ch(A) \: = \: {\mathcal R} \: \tr_s \left( e^{- (A^2 - l)} \right) \in
\Omega^*_{ab}.
\end{equation}
\end{definition}

We have the following analog of Lemma \ref{lemma8}.

\begin{lemma}
Given ${\mathcal E}$, if
$\eta$ is a closed graded trace on $\Omega^*$
then $\eta(\ch(A))$ is independent
of the choice of $A$. If $\eta_1$ and $\eta_2$ are homologous closed
graded traces then $\eta_1(\ch(A)) \: = \: \eta_2(\ch(A))$.
\end{lemma}

\section{Differential Calculus for Foliation
Groupoids}\label{differential}

In this section, given a
foliation groupoid $G$, we construct a graded
algebra $\Omega^{*}{\mathcal B}$ whose degree-$0$ component
${\mathcal B}$ is the convolution algebra of $G$.
We then construct a degree-$1$ derivation
$d = d^H$ of $\Omega^{*}{\mathcal B}$.
Finally, we compute $d^2$.

\subsection{The differential forms}

Let $G$ be a groupoid. We use the groupoid notation of
\cite[Section II.5]{Connes (1994)}.
The units of $G$ are denoted $G^{(0)}$ and the
range and source maps
are denoted $r, s \: : \: G \rightarrow G^{(0)}$.
To construct the product of $g_0, g_1 \in G$, we must have
$s(g_0) \: = \: r(g_1)$. Then $r(g_0 g_1) \: = \: r(g_0)$ and
$s(g_0 g_1) \: = \: s(g_1)$. Given $m \in G^{(0)}$, put
$G^m \: = \: r^{-1}(m)$, $G_m \: = \: s^{-1}(m)$ and $G^m_m \: = \:
G^m \cap G_m$.

We assume that $G$ is a Lie groupoid, meaning that $G$ and $G^{(0)}$ are
smooth manifolds, and $r$ and $s$ are smooth submersions.
For simplicity we will assume that $G$ is Hausdorff.
The
results of the paper extend to the nonHausdorff case, using the
notion of differential forms on a nonHausdorff manifold given by
Crainic and Moerdijk
\cite[Section 2.2.5]{Crainic-Moerdijk (1999)}. (The paper
\cite{Crainic-Moerdijk (1999)} is an extension of work by
Brylinski and Nistor
\cite{Brylinski-Nistor (1994)}.)

The Lie algebroid ${\mathfrak g}$ of $G$ is a vector bundle
over $G^{(0)}$ with fibers ${\mathfrak g}_m \: = \: T_m G^m \: = \:
\Ker(dr_m \: : \: T_m G \rightarrow  T_m G^{(0)})$.   The anchor
map ${\mathfrak g} \rightarrow TG^{(0)}$, a map of vector bundles,
is the restriction of
$ds_m \: : \: T_m G \rightarrow  T_m G^{(0)}$ to ${\mathfrak g}_m$.
In general, the image of the anchor map need not be of constant rank.

We now assume that $G$ is a foliation
groupoid in the sense of
\cite{Crainic-Moerdijk (2001)}, i.e. that $G$ satisfies one of the three
following equivalent conditions \cite[Theorem 1]{Crainic-Moerdijk
(2001)} :
\\ 1. $G$ is Morita equivalent to a smooth \'etale groupoid. \\
2. The anchor map of $G$ is injective. \\
3. All isotropy Lie groups $G_m^m$ of $G$ are discrete. \\ \\
{\bf Example 4 : } If $G$ is an smooth \'etale groupoid then $G$ is a
foliation groupoid.  If $(M, {\mathcal F})$ is a smooth foliated manifold
then its holonomy groupoid 
(see Connes \cite[Section II.8.$\alpha$]{Connes (1994)}) and its
monodromy (= homotopy) groupoid
(see Baum-Connes \cite{Baum-Connes (1985)} and
Phillips \cite{Phillips (1987)}) are foliation
groupoids. In this case, the anchor map is the inclusion map
$T{\mathcal F} \rightarrow TM$.
If a Lie group $L$ acts smoothly on
a manifold $M$ and the isotropy groups 
$L_m \: = \: \{l \in L \: : \: ml \: = \: m\}$ are discrete 
then the cross-product groupoid $M \rtimes L$ is a foliation groupoid.
\\

Put $M \: =\: G^{(0)}$.  It inherits a foliation ${\mathcal F}$, with
the leafwise tangent bundle $T{\mathcal F}$ being the image of the anchor map.

Note that the foliated manifold $(M, {\mathcal F})$ has a holonomy groupoid
${\Hol}$ which is itself a foliation groupoid.
However, ${\Hol}$ may not be the same as $G$.  If $G$ is a
foliation groupoid
with the property that $G_m$ is connected for all $m$ then $G$
lies between the holonomy groupoid of ${\mathcal F}$ and the monodromy
groupoid of ${\mathcal F}$; see \cite[Proposition
1]{Crainic-Moerdijk (2001)} for further discussion. The reader may just
want to keep in mind the case when $G$ is actually the holonomy groupoid
of a foliated manifold $(M, {\mathcal F})$.

Let $\tau \: = \: TM/T{\mathcal F}$ be the normal bundle to the foliation.
Given $g \in G$, let $U \subset M$ be a sufficiently small
neighborhood of $s(g)$ and let $c \: : \: U \rightarrow G$ be a smooth
map such that $c(s(g)) \: = \: g$ and $s \circ c \: = \: \Id_U$.
Then $d(r \circ c)_{s(g)} \: : \: T_{s(g)}M \rightarrow T_{r(g)}M$ sends
$T_{s(g)}{\mathcal F}$ to $T_{r(g)}{\mathcal F}$. The induced map
from $\tau_{s(g)}$ to $\tau_{r(g)}$ has an inverse
$g_* \: : \: \tau_{r(g)} \rightarrow \tau_{s(g)}$ called the holonomy of
the element $g \in G$. It is independent
of the choices of $U$ and $c$.

Let ${\mathcal D}$  denote the real line bundle on $M$
formed by leafwise densities.
We define a graded algebra $\Omega^*{\mathcal B}$
whose components, as vector spaces, are given by
\begin{equation}\label{3.4}
\Omega^{n}{\mathcal B} \: = \: C^{\infty}_c\left(
G; \Lambda^n (r^* \tau^*)
\otimes s^*{\mathcal D} \right)
\end{equation}
In particular,
\begin{equation}\label{3.5}
{\mathcal B} \: = \: \Omega^{0}{\mathcal B} \: = \:
C^\infty_c(G; s^* {\mathcal D})
\end{equation}
is the groupoid algebra.
(Instead of using half-densities, we have placed a full density
at the source.)
The product of $\phi_1 \in \Omega^{n_1}{\mathcal B}$ and
$\phi_2 \in \Omega^{n_2}{\mathcal B}$ is given by
\begin{equation}
(\phi_1 \phi_2)(g) \: = \:
 \int \limits_{g^\prime g^{\prime \prime} \: = \: g}
\phi_1(g^\prime) \: \wedge \:
\phi_2(g^{\prime \prime}).
\end{equation}
In forming the wedge product, the holonomy of $g^\prime$ is  used
to identify conormal spaces.

Let $T^H M$ be a
horizontal distribution on $M$, i.e. a splitting of
the short exact sequence $0 \rightarrow T{\mathcal F} \rightarrow TM
\rightarrow \tau \rightarrow 0$.
Then there is a horizontal differentiation  $d^H \: : \:
\Omega^{n}{\mathcal B} \rightarrow
\Omega^{n+1}{\mathcal B}$, which we now define. The
definition will proceed by building up
$d^H$ from smaller pieces
(compare \cite[Section II.7.$\alpha$, Proposition 3]{Connes (1994)}).

First, the choice of horizontal distribution allows us to define a
horizontal differential $d^H \: : \: \Omega^*(M) \rightarrow
\Omega^{*+1}(M)$  as in 
Bismut-Lott \cite[Definition 3.2]{Bismut-Lott (1995)}
and Connes
\cite[Section III.7.$\alpha$]{Connes (1994)}.
Using the local description
of an element of
$C^{\infty}\left(M; {\mathcal
D}\right)$ as a vertical $\dim({\mathcal F})$-form on $M$, we
also obtain a horizontal differential $d^H \: : \:
C^{\infty}\left(M; {\mathcal D}\right)
\rightarrow
C^{\infty}\left(M; \tau^*\otimes {\mathcal D}\right)$
\cite[Section III.7.$\alpha$]{Connes (1994)} and a
horizontal differential
$d^H \: : \: C^{\infty}\left(M; \Lambda^n \tau^*\right) \rightarrow
C^{\infty}\left(M; \Lambda^{n+1} \tau^*\right)$.

Given $f \in C^\infty_c \left( G \right)$, we
now define its horizontal differential $d^H f \in
C^\infty_c \left( G; r^* \tau^* \right)$ by simultaneously
differentiating $f$ with respect to its arguments, in a horizontal
direction. That is, consider a
point $g \in
G$ and a vector $X_0\in \tau_{r(g)}$.
Put $X_1 \: = \: g_*(X_0)$.
Next, use the horizontal distribution $T^H M$ to construct
the corresponding horizontal vectors $\widetilde{X}_0$ and
$\widetilde{X}_1$. We now have a vector
$\widetilde{X} \: = \: \left(\widetilde{X}_0,
\widetilde{X}_1\right)\in T_{(r(g), s(g))}(M \times M)$. It
is the image of a unique vector $X \in T_{g}G$
under the immersion
\begin{equation} \label{3.18}
(r, s) \: : \: G
\to M \times M.
\end{equation}
We define $d^H f$ by putting
$\left( (d^H f)(X_0) \right)(g) \: = \: Xf$.

Next, to horizontally differentiate an element of
$C^\infty_c \left( G; \Lambda^n(r^* \tau^*)
\otimes s^* {\mathcal D} \right)$,
we write it as  a finite sum of terms of the form $f
\: r^*(\omega) \: s^*(\beta)$,
with $f \in C^\infty_c \left( G \right)$,
$\omega
\in C^\infty (M; \Lambda^n
\tau^*)$,
and
$\beta
\in C^\infty(M; {\mathcal D})$.
For an element of this form, put
\begin{equation}
d^H \left(f \: r^*(\omega) \: s^*(\beta) \right)
 \: = \:
(d^H f) \: r^*(\omega) \: s^*(\beta) \:  + \:
f \: r^*(d^H\omega) \: s^*(\beta) \: + \:
(-1)^n \: f \: r^*(\omega) \: s^*(d^H \beta),
\end{equation}
where the holonomy is used in defining products.

\begin{lemma}
The operator $d^H$ is a graded derivation of $\Omega^* {\mathcal B}$.
\end{lemma}
\begin{pf}
This follows from a straightforward computation, which we omit.
\end{pf}

Put
$d \: = \: d^H$.
We now describe $\alpha \: = \: d^2$. Let $T \in \Omega^2(M; T{\mathcal F})$
be the curvature of the horizontal distribution $T^HM$
\cite[(3.11)]{Bismut-Lott (1995)}.
It is a
horizontal $2$-form on $M$ with values in $T {\mathcal F}$, defined by
$T(X_1, X_2) \: = \: - \: P^{vert} \: [X_1^H, X_2^H]$.
One can define the Lie derivative
${\mathcal L}_T \: : \: \Omega^*(M) \rightarrow \Omega^{*+2}(M)$, an operation
which increases the horizontal grading by two, as in
\cite[(3.14)]{Bismut-Lott (1995)}. Then one can define
${\mathcal L}_T \: : \: C^\infty(M; {\mathcal D}) \rightarrow
C^\infty(M; \Lambda^2 \tau^* \otimes {\mathcal D})$ and
${\mathcal L}_T \: : \: C^\infty(M; \Lambda^n \tau^*) \rightarrow
C^\infty(M; \Lambda^{n+2} \tau^*)$ in obvious ways.

Given $f \in C^\infty_c \left( G \right)$, we define its
Lie derivative ${\mathcal L}_T f \in C^\infty_c \left( G;
\Lambda^2 (r^* \tau^*) \right)$ by simultaneously
differentiating $f$ with respect to its arguments, in the
vertical direction. That is, consider a point $g
\in G$ and $X_0, Y_0 \in \tau_{r(g)}$.
Put $X_1 \: = \: g_*(X_0)$ and $Y_1 \: = \: g_*(Y_0)$.
Next, use the horizontal distribution $T^H M$
to construct the corresponding horizontal vectors
$\widetilde{X}_0$, $\widetilde{X}_1$, $\widetilde{Y}_0$
and $\widetilde{Y}_1$. Consider the vertical
vectors $T(\widetilde{X}_0, \widetilde{Y}_0) \in T_{r(g)}
{\mathcal F}$ and
$T(\widetilde{X}_1, \widetilde{Y}_1) \in T_{s(g)}
{\mathcal F}$.
We now have a total vector $\widetilde{V} \: = \:
\left( T(\widetilde{X}_0, \widetilde{Y}_0),
T(\widetilde{X}_1, \widetilde{Y}_1) \right) \in T_{(r(g),s(g))}
(M \times M)$. It is the image of a unique
vector $V \in T_{g}G$ under the
immersion (\ref{3.18}). We define ${\mathcal L}_T f$ by putting
$\left( ({\mathcal L}_T f)(X_0, Y_0)  \right) (g)\:
= \: Vf$.

Now for $f \: r^*(\omega) \: s^*(\beta)$ as before, we put
\begin{equation}
{\mathcal L}_T \left(f \: r^*(\omega) \: s^*(\beta)
\right) \: = \:
({\mathcal L}_T f) \: r^*(\omega) \: s^*(\beta)
\: + \:
f \: r^*({\mathcal L}_T \omega) \: s^*(\beta)
\: + \:
f \: r^*(\omega) \:
s^*({\mathcal L}_T \alpha_{1}),
\end{equation}
where the holonomy is used in defining products.

\begin{lemma}
We have
\begin{equation}
\alpha \: = \: - \: {\mathcal L}_T.
\end{equation}
\end{lemma}
\begin{pf}
This follows from
the method of proof of \cite[(3.13)]{Bismut-Lott (1995)} or
\cite[Section III.7.$\alpha$]{Connes (1994)}.
\end{pf}
\noindent
{\bf Remark : } One
can consider $\alpha$ to be commutation with a (distributional)
element of the multiplier algebra $C^{- \infty} \left (G;
\Lambda^2 (p_0^* \tau^*) \otimes p_1^*{\mathcal D} \right)$,
namely the one that implements the Lie differentiation
\cite[Section III.7.$\alpha$, Lemma 4]{Connes (1994)}.

\section{Superconnection and  Chern
character}\label{superconnection}

In this section we consider a smooth manifold $P$ on which $G$ acts
freely, properly and cocompactly, along with a $G$-invariant
$\Z_2$-graded vector bundle $E$ on $P$.  We construct a corresponding
left-${\mathcal B}$-module ${\mathcal E}$.
Given a $G$-invariant
Dirac-type operator which acts on sections of $E$, we consider the
Bismut
superconnections $\{A_s\}_{s > 0}$. We compute the $s \rightarrow 0$
limit of the pairing between the Chern character of $A_s$ and
a closed graded trace on $\Omega^* {\mathcal B}$ that is concentrated on
the units $M$. More detailed
summaries appear at the beginnings of the subsections.

\subsection{Module and Connection}\label{subsection5.1}

In this subsection we consider a left ${\mathcal B}$-module
${\mathcal E}$ consisting of sections of $E$, and its extension to a
left $\Omega^*{\mathcal B}$-module $\Omega^*{\mathcal E}$. We construct
a map $l : {\mathcal E} \rightarrow \Omega^2{\mathcal E}$
satisfying (\ref{l}). Given a lift
$T^HP$ of $T^HM$, we construct a connection
$\nabla^{\mathcal E}$ on ${\mathcal E}$.

Let $P$ be a smooth $G$-manifold \cite[Section II.10.$\alpha$,
Definition 1]{Connes (1994)}. That is, first of all, there is a submersion
$\pi \: : \: P \rightarrow M$. Given $m \in M$, we write
$Z_m \: = \: \pi^{-1}(m)$. Putting
\begin{equation}
P \times_r G \: = \: \{ (p, g) \in P \times G \: : \:
p \: \in \: Z_{r(g)} \},
\end{equation}
we must also have a smooth map $P \times_r G \rightarrow P$, denoted
$(p, g) \rightarrow p g$, such that
$p  g \: \in \: Z_{s(g)}$ and
$(p g_1) g_2 \: = \:
p (g_1 g_2)$ for all $(g_1, g_2) \in G^{(2)}$.
It follows that for each $g \in G$, the map $p \rightarrow p  g$
gives a diffeomorphism from $Z_{r(g)}$ to
$Z_{s(g)}$. Let
${\mathcal D}_Z$ denote the real line bundle on $P$
formed by the fiberwise densities.

Hereafter we assume that $P$ is a proper $G$-manifold
\cite[Section II.10.$\alpha$, Definition 2]{Connes (1994)},
i.e. that the map $P \times_{r} G
\rightarrow P \times P$ given by $(p, g) \rightarrow (p, p g)$
is proper.
We also assume that $G$ acts cocompactly on $P$, i.e. that
the quotient of $P$ by the equivalence relation ($p \sim p^\prime$ if
$p \: = \: p^\prime g$ for some $g \in G$) is compact.
And we assume that $G$ acts freely on
$P$, i.e. that $pg \: = \: p$ implies that $g \in M$.
Then $P/G$ is a smooth compact manifold.
\\ \\
{\bf Example 5 : } Take $P \: = \: G$, with $\pi \: = \: s$. Then
$G$ acts properly, freely, and, if $M$ is compact, cocompactly on $P$. \\

We will say that a covariant object (vector bundle, connection,
metric, etc.) on $P$ is $G$-invariant if it is the pullback of a similar
object from $P/G$.
Let $E$ be a $G$-invariant $\Z_2$-graded vector bundle on $P$, with
supertrace $\tr_s$ on $\End(E)$.
Put
${\mathcal E} = C^{\infty}_c(P; E)$.
It is a left-${\mathcal B}$-module, with the action of $b \in {\mathcal B}$ on
$\xi \in {\mathcal E}$ given by
\begin{equation}\label{5.1}
(b \xi)(p) \: = \: \int_{G^{\pi(p)}} b(g) \: \xi(pg).
\end{equation}
In writing (\ref{5.1}), we have used the $g$-action to identify
$E_p$ and $E_{pg}$.

Put
\begin{equation}
\Omega^{n}{\mathcal E} = C^{\infty}_c \left( P ;
\Lambda^n
(\pi^* \tau^*) \otimes  E
\right).
\end{equation}
Then $\Omega^{*}{\mathcal E}$ is a left-$\Omega^*{\mathcal B}$-module with the action
of $\Omega^{*}{\mathcal B}$ on $\Omega^{*}{\mathcal E}$
given by
\begin{equation}
(\phi  \: \omega)(p) \: = \:
\int_{G^{\pi(p)}} \phi(g)\: \wedge \: \omega(pg).
\end{equation}

Let $\widetilde{\mathcal F}$ be the foliation on $P$ whose leaf through
$p \in P$ consists of the elements $pg$ where $g$ runs through the
connected component of $G^{\pi(p)}$ that contains the unit $\pi(p)$.
Note that $\dim(\widetilde{\mathcal F}) \: = \: \dim({\mathcal F})$.
Given $p \in P$ and $X, Y \in \tau_{\pi(p)}$, let
$\widetilde{T}({X}, {Y}) \in T_p \widetilde{\mathcal F}$ be the lift of
${T}({X}, {Y}) \in T_{\pi(p)} {\mathcal F}$.
Define $l \: : \: {\mathcal E} \rightarrow
\Omega^{2}{\mathcal E}$ by saying that for $X, Y \in \tau_{\pi(p)}$
and $\xi \in {\mathcal E}$,
\begin{equation} \label{5.5}
(l(\xi)(X, Y))(p) \: = \: - \: \widetilde{T}
({X}, {Y}) \xi.
\end{equation}
Here we have used the $G$-invariance of $E$
to define the action of
$\widetilde{T}({X}, {Y})$ on $\xi$.
\begin{lemma} \label{lemma16}
For all $X, Y \in \tau_{\pi(p)}$,
$b \in {\mathcal B}$ and $\xi \in {\mathcal E}$,
\begin{equation} \label{5.6}
l(b \xi) \: = \: \alpha(b) \: \xi \: + \: b \: l(\xi).
\end{equation}
\end{lemma}
\begin{pf}
We have
\begin{equation}
(l(b \xi)(X, Y))(p) \: = \: - \: \widetilde{T}(X,Y)
\int_{G^{\pi(p)}} b(g) \: \xi(pg) \: = \: - \:
\int_{G^{\pi(p)}}  {T}(X,Y)b(g) \: \xi(pg),
\end{equation}
\begin{equation}
(\alpha(b)(X, Y) \xi)(p) \: = \: - \:
\int_{G^{\pi(p)}} \left( T(X, Y) b \: + \: T(g_*X, g_*Y) b \right)(g) \:
\xi(pg)
\end{equation}
and
\begin{equation}
(b l(\xi)(X, Y))(p) \: = \: - \:
\int_{G^{\pi(p)}} b(g) \: \widetilde{T}(g_*X, g_*Y) \xi(pg).
\end{equation}
Then
\begin{align}\label{5.10}
& (l(b \xi)(X, Y))(p) \: - \: (\alpha(b)(X, Y) \xi)(p) \: - \:
(b l(\xi)(X, Y))(p) \: = \\
& \int_{G^{\pi(p)}} \left( T(g_*X, g_*Y) b(g) \:
\xi(pg) \: + \: b(g) \: \widetilde{T}(g_*X, g_*Y) \xi(pg) \right). \notag
\end{align}
We can write (\ref{5.10}) more succinctly as
\begin{equation}\label{5.11}
l(b \xi) \: - \: \alpha(b) \xi \: - \:
b l(\xi) \: =
\int_{G^{\pi(p)}} {\mathcal L}_{\widetilde{T}} (b(g) \xi(pg)),
\end{equation}
where the Lie differentiation is at $pg$. The right-hand-side of
(\ref{5.11}) vanishes, being the integral of a Lie derivative of a
compactly-supported density.
\end{pf}

We extend $l$ to a linear map
$l \: : \: \Omega^{n} {\mathcal E} \rightarrow
\Omega^{n+2} {\mathcal E}$ as Lie differentiation
in the $\widetilde{T}$-direction with respect to $P$.
\begin{lemma} \label{lll}
For all $\omega \in \Omega^* {\mathcal B}$ and
$\mu \in \Omega^* {\mathcal E}$,
\begin{equation}
l(\omega \mu) \: = \: \alpha(\omega) \: \mu \: + \: \omega \: l(\mu).
\end{equation}
\end{lemma}
\begin{pf}
The proof is similar to that of Lemma \ref{lemma16}. We omit the details.
\end{pf}

There is a pullback foliation $\pi^* {\mathcal F}$ on $P$ with the
same codimension as ${\mathcal F}$,  satisfying
$T\pi^* {\mathcal F} \: = \: (d\pi)^{-1} T{\mathcal F}$.
Let $\mu \: : \: P \rightarrow
P/G$ be the quotient map.
Then $P/G$
is a smooth compact manifold with a foliation $F \: = \: (\pi^*
{\mathcal F})/G$ satisfying
$(d\mu)^{-1} TF \: =  \: T \pi^*
{\mathcal F}$. We note that the normal bundle $NF$ to $F$ satisfies
$\mu^* NF \: = \: \pi^* \tau$.

Let $T^H (P/G)$ be a horizontal distribution on $P/G$, transverse to
$F$. Then $(d\mu)^{-1}(T^H (P/G))$ is a $G$-invariant
distribution on $P$ that is transverse
to the vertical tangent bundle $TZ$.
Put $T^H P \: = \: (d\mu)^{-1}(T^H (P/G)) \cap (d\pi)^{-1} (T^H M)$, a
distribution
on $P$ that is transverse to $\pi^* {\mathcal F}$ and that projects
isomorphically
under $\pi$ to $T^HM$.

Let
$\nabla^{\mathcal E} : {\mathcal E} \to
\Omega^{1} {\mathcal E}$
be covariant differentiation on ${\mathcal E} \: = \:
C^\infty_c(P; E)$ with respect
to $T^H P$.

\begin{lemma} \label{llll}
$\nabla^{\mathcal E}$ is a connection.
\end{lemma}
\begin{pf}
We wish to show that
\begin{equation} \label{wanttoshow}
\nabla^{\mathcal E}(b \xi)=b \nabla^{\mathcal E} \xi +(d^H b) \xi.
\end{equation}
As the claim of the lemma is local on $P$, consider first the case
when $T^H(P/G)$ is integrable. Let $T^H P_1$ and $\nabla^{\mathcal E}_1$
denote the corresponding objects on $P$. Then one is geometrically in a product
situation and one can reduce to the case $P \: = \: M$,
where one can check that (\ref{wanttoshow}) holds.
If $T^H(P/G)$ is not integrable then
$T^HP \: - \: T^H P_1 \in \Hom(\pi^* \tau, TZ)$ is the
pullback under $\mu$ of an element of $\Hom(NF, TF)$. Hence
$T^HP \: - \: T^H P_1$ is $G$-invariant and it follows that
$\nabla^{\mathcal E} \: - \: \nabla^{\mathcal E}_1$ commutes with
${\mathcal B}$, which proves the lemma.
\end{pf}

We extend
$\nabla^{\mathcal E}$ to act on
$\Omega^{*} {\mathcal E}$
so as to satisfy Leibnitz' rule.

\begin{lemma}
For all $\xi \in {\mathcal E}$,
\begin{equation}
l(\nabla^{\mathcal E} \xi) \: = \: \nabla^{\mathcal E} l(\xi).
\end{equation}
\end{lemma}
\begin{pf}
As $d^H$ commutes with $(d^H)^2$, it follows that
$d^H$ commutes with ${\mathcal L}_T$. As the claim of the lemma is
local on $P$, consider first the case when $T^H (P/G)$ is integrable.
Let $T^H P_1$ and $\nabla^{\mathcal E}_1$ denote the corresponding
objects on $P$.
Then one is in a local product situation and the lemma follows from the
fact that $d^H$ commutes with ${\mathcal L}_T$. If
$T^H (P/G)$ is not integrable then
$\nabla^{\mathcal E} \: - \: \nabla^{\mathcal E}_1$ is given by covariant
differentiation in the $TZ$ direction,
with respect to $T^HP \: - \: T^H P_1 \in \Hom(\pi^* \tau, TZ)$.
As
$\widetilde{T}$ pulls back from $M$,
$\nabla^{\mathcal E} \: - \: \nabla^{\mathcal E}_1$ commutes with $l$.
The lemma follows.
\end{pf}

\subsection{Supertraces}\label{subsection5.2}

In this subsection we consider a certain
algebra $\End^\infty_{\mathcal B} ({\mathcal E})$ of operators with smooth
kernel on $P$. We show that a trace on ${\mathcal B}$, concentrated on the
units $M$, gives a supertrace
on $\End^\infty_{\mathcal B} ({\mathcal E})$. We then consider an algebra
$\Hom^\infty_{\mathcal B} ({\mathcal E}, \Omega^* {\mathcal E})$ of
form-valued operators. We show that a closed graded trace on
$\Omega^*{\mathcal B}$, concentrated on $M$,
gives rise to a closed graded trace on
$\Hom^\infty_{\mathcal B} ({\mathcal E}, \Omega^* {\mathcal E})$.

An operator $K \in \End_{\mathcal B}( {\mathcal E})$ has a Schwartz kernel
$K(p^\prime | p)$ so that
\begin{equation}\label{5.20}
(K\xi)(p) \: = \: \int_{Z_{\pi(p)}} \xi(p^\prime) \: K(p^\prime | p).
\end{equation}
Define $q^\prime, q \: : \: P \times_{M} P \rightarrow P$ by
$q^\prime(p^\prime, p) \: = \: p^\prime$ and $q(p^\prime, p) \: = \: p$.
Let $\End^\infty_{\mathcal B} ({\mathcal E})$ denote the subalgebra of
$\End_{\mathcal B}( {\mathcal E})$ consisting of
operators whose Schwartz kernel lies in
$C^\infty_c( P \times_{M} P; (q^\prime)^* {\mathcal D}_Z \otimes
\Hom((q^\prime)^* E, q^* E))$.

Choose $\Phi \in
C^\infty_c(P; \pi^*{\mathcal D})$ so that
\begin{equation}\label{5.21}
\int_{G^{\pi(p)}} \Phi(pg) = 1
\end{equation}
for all $p\in P$; that such a $\Phi$ exists was shown by Tu 
\cite[Proposition 6.11]{Tu (1999)}.
Define $\tau K \in C^\infty_c(M; {\mathcal D})$ by
\begin{equation}
(\tau K)(m) \: = \: \int_{Z_m} \Phi(p) \: \tr_s K(p|p).
\end{equation}

\begin{proposition}\label{prop3}
Let $\rho$ be a linear functional on
$C^\infty_c(M; {\mathcal D})$. Suppose that
the linear functional
$\eta$ on ${\mathcal B}$, defined by
\begin{equation}
\eta(b) \: = \: \rho( b \big|_{M}),
\end{equation}
is a trace on ${\mathcal B}$.
Then $\rho \circ \tau$ is a supertrace on
$\End^\infty_{\mathcal B} ({\mathcal E})$.
\end{proposition}
\begin{pf}
Consider the algebra $\End_{C^\infty_c(M)}( {\mathcal E})$.
An operator $K \in
\End_{C^\infty_c(M)}( {\mathcal E})$ has a Schwartz kernel
$K(p | p^\prime)$ so that
\begin{equation}
(K \xi)(p) \: = \: \int_{Z_{\pi(p)}} K(p|p^\prime) \: \xi(p^\prime).
\end{equation}
(Note the difference in ordering as compared to (\ref{5.20}).) For
this proof, define $q, q^\prime \: : \: P \times_{M} P \rightarrow
P$ by $q(p, p^\prime) \: = \: p$ and $q^\prime(p, p^\prime) \:
= \: p^\prime$. Let $\End^\infty_{C^\infty_c(M)} ({\mathcal E})$
denote the subalgebra of $\End_{C^\infty_c(M)}( {\mathcal E})$
consisting of operators whose Schwartz kernel lies in $C^\infty_c(
P \times_{M} P; q^* {\mathcal D}_Z \otimes \Hom((q^\prime)^* E,
q^* E)$. The product in $\End^\infty_{C^\infty_c(M)}
({\mathcal E})$ is given by
\begin{equation}
(K K^\prime)(p | p^\prime) \: = \: \int_{p^{\prime \prime}}
K(p | p^{\prime \prime}) \:
K^\prime(p^{\prime \prime} | p^\prime).
\end{equation}
Note that an element of
$\End^\infty_{C^\infty_c(M)} ({\mathcal E})$ is not necessarily
$G$-invariant.
Note also that there is an
injective homomorphism
$\End^\infty_{\mathcal B} ({\mathcal E}) \rightarrow
\End^\infty_{C^\infty_c(M)} ({\mathcal E})^{op}$, where ${op}$ denotes
the opposite algebra, i.e. with the transpose multiplication.
There is a fiberwise $G$-invariant supertrace
$Tr_s \: : \:
\End^\infty_{C^\infty_c(M)} ({\mathcal E}) \rightarrow C^\infty_c(M)$ given by
\begin{equation}\label{5.26}
(Tr_s \: K)(m) \: = \: \int_{Z_m} \tr_s K(p|p).
\end{equation}

Consider  the algebra ${\mathcal B}
\otimes_{C^{\infty}_c(M)} \End^\infty_{C^\infty_c(M)} ({\mathcal E})$.
The product in the algebra takes into
account the action of ${\mathcal B}$ on $\End^\infty_{C^\infty_c(M)} ({\mathcal E})$,
which derives from the $G$-action on $P$.
An element of the
algebra has a kernel
$K(g, p | p^\prime)$, where $p, p^\prime \in
Z_{s(g)}$. The product is given by
\begin{equation}
(K_1 K_2)(g, p | p^\prime) \: = \:
\int_{g^\prime g^{\prime \prime} \: = \: g}
\int_{p^{\prime \prime} \in Z_{s(g^\prime)}}
K_1(g^\prime, p (g^{\prime \prime})^{-1} | p^{\prime \prime}) \:
K_2(g^{\prime \prime}, p^{\prime \prime} g^{\prime \prime} |
p^{\prime}).
\end{equation}

The supertrace (\ref{5.26}) induces a map $Tr_s \: : \: {\mathcal
B} \otimes_{C^{\infty}_c(M)} \End^\infty_{C^\infty_c(M)}
({\mathcal E}) \rightarrow {\mathcal B}$ by
\begin{equation}
(Tr_s \: K)(g) \: = \: \int_{Z_{s(g)}} \tr_s K(g, p|p).
\end{equation}
\begin{lemma}\label{lemma19}
$\eta \circ Tr_s$ is a supertrace on
${\mathcal B}
\otimes_{C^{\infty}_c(M)} \End^\infty_{C^\infty_c(M)} ({\mathcal E})$.
\end{lemma}
\begin{pf}
We can formally write
\begin{equation}
(\eta \circ Tr_s)(K) \: = \: \int_M \rho(m) \int_{Z_m} \: \tr_s K(m, p|p),
\end{equation}
keeping in mind that $\rho$ is actually distributional.
Then
\begin{align}
(\eta \circ Tr_s)(K_1 K_2) \: & = \:
\int_{g^\prime \in G}
\int_{p \in Z_{r(g^\prime)}}
\int_{p^{\prime \prime} \in Z_{s(g^\prime)}}
\rho(r(g^\prime)) \:
\tr_s \left( K_1(g^\prime, p g^{\prime} | p^{\prime \prime}) \:
K_2((g^{\prime})^{-1}, p^{\prime \prime} (g^{\prime})^{-1} | p) \right) \\
& = \:
\int_{g^\prime \in G}
\int_{p \in Z_{r(g^\prime)}}
\int_{p^{\prime \prime} \in Z_{s(g^\prime)}}
\rho(r(g^\prime)) \:
\tr_s \left(
K_2((g^{\prime})^{-1}, p^{\prime \prime} (g^{\prime})^{-1} | p) \:
K_1(g^\prime, p g^{\prime} | p^{\prime \prime}) \right) \notag \\
& = \:
\int_{g^\prime \in G}
\int_{p^{\prime \prime} \in Z_{r(g^\prime)}}
\int_{p \in Z_{s(g^\prime)}}
\rho(s(g^\prime)) \: \tr_s \left(
K_2(g^{\prime}, p^{\prime \prime} g^{\prime} | p) \:
K_1((g^\prime)^{-1}, p (g^{\prime})^{-1} | p^{\prime \prime}) \right). \notag
\end{align}
However, the fact that $\eta$ is a trace on ${\mathcal B}$
translates into the fact that
\begin{equation}
\int_{g \in G} \rho(s(g)) \: f(g) \: = \:
\int_{g \in G} \rho(r(g)) \: f(g)
\end{equation}
for all $f \in C^\infty_c(G)$, from which the lemma follows.
\end{pf}

We define a map
$i \: : \: \End^\infty_{\mathcal B} ({\mathcal E}) \rightarrow
\left( {\mathcal B}
\otimes_{C^{\infty}_c(M)} \End^\infty_{C^\infty_c(M)}
({\mathcal E})\right)^{op}$ by
\begin{equation}
(i(K))(g, p|p^\prime) \: = \: \Phi(p g^{-1} )K(p | p^\prime).
\end{equation}
\begin{lemma}\label{lemma20}
 The map $i$ is a homomorphism.
 \end{lemma}
\begin{pf}
Given $K_1, K_2 \in \End^\infty_{\mathcal B} ({\mathcal E})$, we have
\begin{align}
\left( i(K_1) \: i(K_2) \right)(g, p | p^\prime) \: & = \:
\int_{g^\prime g^{\prime \prime} \: = \: g} \int_{Z_{s(g^\prime)}}
i(K_1)(g^\prime, p (g^{\prime \prime})^{-1} | p^{\prime \prime}) \:
i(K_2)(g^{\prime \prime}, p^{\prime \prime} g^{\prime \prime} |
p^{\prime}) \\
& = \:
\int_{g^\prime g^{\prime \prime} \: = \: g} \int_{Z_{s(g^\prime)}}
\Phi(pg^{-1}) \:
K_1(p (g^{\prime \prime})^{-1} | p^{\prime \prime}) \:
\Phi(p^{\prime \prime}) \:
K_2(p^{\prime \prime} g^{\prime \prime} |
p^{\prime}) \notag \\
& = \:
\int_{g^\prime g^{\prime \prime} \: = \: g} \int_{Z_{s(g^\prime)}}
\Phi(pg^{-1}) \:
K_1(p | p^{\prime \prime} g^{\prime \prime}) \:
\Phi(p^{\prime \prime}) \:
K_2(p^{\prime \prime} g^{\prime \prime} |
p^{\prime}) \notag \\
& = \:
\int_{g^\prime g^{\prime \prime} \: = \: g}
\int_{Z_{s(g^{\prime \prime})}}
\Phi(pg^{-1}) \:
K_1(p | p^{\prime \prime}) \:
\Phi(p^{\prime \prime} ( g^{\prime \prime})^{-1}) \:
K_2(p^{\prime \prime} |
p^{\prime}) \notag \\
& = \:
\Phi(pg^{-1}) \: \int_{Z_{s(g)}}
K_1(p | p^{\prime \prime}) \:
K_2(p^{\prime \prime} |
p^{\prime}) \notag \\
& = \:
(i(K_2 K_1))(g, p | p^\prime). \notag
\end{align}
Thus $i$ gives a homomorphism from
$\End^\infty_{\mathcal B} ({\mathcal E})^{op}$ to
${\mathcal B}
\otimes_{C^{\infty}_c(M)} \End^\infty_{C^\infty_c(M)}
({\mathcal E})$, from which the lemma follows.
\end{pf}

\begin{lemma}\label{lemma21}
We have
$\eta \circ Tr_s \circ i \: = \: \rho \circ \tau$.
\end{lemma}
\begin{pf}
Given $K \in \End^\infty_{\mathcal B} ({\mathcal E})$, we have
\begin{align}
(\eta \circ Tr_s \circ i)(K) \: & = \:
\int_M \rho(m) \int_{Z_m} \tr_s (i(K))(m, p | p) \\
& = \:
\int_M \rho(m) \int_{Z_m} \Phi(p) \: \tr_s K(p|p) \: = \:
(\rho \circ \tau)(K). \notag
\end{align}
This proves the lemma.
\end{pf}

Proposition \ref{prop3} now follows from Lemmas
\ref{lemma19}-\ref{lemma21}.
\end{pf}

\noindent
{\bf Example 6 : } Let $\mu$ be a holonomy-invariant
transverse measure for ${\mathcal F}$.
Let $\{U_i\}_{i=1}^N$ be an open covering of $M$
by flowboxes, with $U_i \: = \: V_i \times W_i$, $V_i \subset
\R^{codim({\mathcal F})}$ and $W_i \subset
\R^{dim({\mathcal F})}$. Let $\mu_i$ be the measure on $V_i$ which
is the restriction of $\mu$. Let $\{ \phi_i \}_{i=1}^N$ be a
partition of unity that is subordinate to $\{U_i\}_{i=1}^N$.
For $f \in C^\infty_c(M; {\mathcal D})$, put $\rho(f) \: = \:
\sum_{i=1}^N \int_{V_i} \left( \int_{W_i} \phi_i \: f \right) \: d\mu_i$.
Then $\rho$ satisfies the hypotheses of Proposition \ref{prop3}. \\

An operator $K \in \Hom_{\mathcal B} ({\mathcal E}, \Omega^{*} {\mathcal E})$
has a Schwartz kernel
$K(p^\prime | p)$ so that
\begin{equation} \label{5.35}
(K\xi)(p) \: = \:
\int_{Z_{\pi(p)}} \xi(p^\prime) \: K(p^\prime |  p).
\end{equation}
Let $\Hom^\infty_{\mathcal B} ({\mathcal E}, \Omega^{n} {\mathcal E})$
denote the subspace of $\Hom_{\mathcal B} ({\mathcal E}, \Omega^{n} {\mathcal E})$
consisting of
operators whose Schwartz kernel lies in
\begin{equation}
C^\infty_c( P \times_M P ;
\Lambda^n((\pi \circ q)^* \tau^*) \otimes (q^\prime)^* {\mathcal D}_Z \otimes
\Hom((q^\prime)^* E, q^* E)).
\end{equation}
Define
$
\tau K \in C^\infty_c(M; \Lambda^n \tau^* \otimes {\mathcal D})$ by
\begin{equation}\label{5.38}
(\tau K)(m) \: = \: \int_{Z_m} \Phi(p) \: \tr_s K(p | p).
\end{equation}

\begin{proposition}\label{prop4}
Let $\rho$ be a linear functional on
$C^\infty_c(M; \Lambda^n \tau^* \otimes {\mathcal D})$
Suppose that the linear functional
$\eta$ on $\Omega^n{\mathcal B}$, defined by
\begin{equation}
\eta(\phi) \: = \: \rho( \phi \big|_M),
\end{equation}
is a graded trace on $\Omega^*{\mathcal B}$.
Then $\rho \circ \tau$
is a graded trace on
$\Hom^\infty_{\mathcal B} ({\mathcal E},
\Omega^{*} {\mathcal E})$.
\end{proposition}
\begin{pf}
The proof is similar to that of
Proposition \ref{prop3}. We omit the details.
\end{pf}

\begin{proposition}\label{prop5}
Let $\rho$ be a linear functional on
$C^\infty_c(M; \Lambda^n \tau^* \otimes {\mathcal D})$
Suppose that the linear functional
$\eta$ on $\Omega^n{\mathcal B}$, defined by
\begin{equation}
\eta(\phi) \: = \: \rho( \phi \big|_M),
\end{equation}
is a closed graded trace on $\Omega^*{\mathcal B}$.
Then $\rho \circ \tau$
annihilates $[\nabla, K]$ for all $K \in
\Hom^\infty_{\mathcal B} ({\mathcal E},
\Omega^{n-1} {\mathcal E})$.
\end{proposition}
\begin{pf}
It suffices to show that
\begin{equation} \label{suffices}
(\rho \circ \tau)
([\nabla^{\mathcal E}, K]) \: = \:
\eta \left( d^H
(\tau (K)) \right).
\end{equation}

Let $\nabla^{{\mathcal E_0}} \: : \: C^\infty_c(P) \rightarrow
C^\infty_c(P; \pi^* \tau^*)$ be differentiation in the $T^HP$-direction.
It follows from (\ref{5.38}) that
\begin{align}
(d^H (\tau K))(m) \: = \: &
 \: \int_{Z_m} \Phi(p) \: \tr_s [\nabla^{\mathcal E}, K]
(p | p) \: + \: \\
& \: \int_{Z_m} \nabla^{{\mathcal E}_0} \Phi(p) \: \wedge \:
\tr_s K(p | p). \notag
\end{align}
Now $\eta \left( \int_{Z_m} \nabla^{{\mathcal E_0}} \Phi(p) \:
\wedge \: \tr_s K(p | p) \right)$
can be written as $\int_P \nabla^{{\mathcal E_0}} \Phi \wedge {\mathcal O}$
for some $G$-invariant ${\mathcal O}$. From (\ref{5.21}),
$\int_{G^{\pi(p)}} \nabla^{{\mathcal E_0}} \Phi(pg) = 0$. Then decomposing
the measure on $P$ with respect to $P \rightarrow P/G$ gives that
$\int_P \nabla^{{\mathcal E_0}} \Phi \wedge {\mathcal O} \: = \: 0$.
Equation (\ref{suffices}) follows.
\end{pf}
\noindent
{\bf Example 7 : } Following the notation of Example 6,
let $c$ be a closed holonomy-invariant
transverse $n$-current for ${\mathcal F}$.
Let $c_i$ be the $n$-current on $V_i$ which
is the restriction of $c$. Let $\{ \phi_i \}_{i=1}^N$ be a
partition of unity that is subordinate to $\{U_i\}_{i=1}^N$.
For
$\omega \in C^\infty_c(M; \Lambda^n \tau^* \otimes {\mathcal D})$,
put $\rho(\omega) \: = \:
\sum_{i=1}^N \langle \left( \int_{W_i} \phi_i \: \omega \right), \: c_i
\rangle$.
Then $\rho$ satisfies the hypotheses of Proposition \ref{prop5}. \\

\subsection{The $s \rightarrow 0$ limit of the
superconnection Chern character}\label{subsection5.3}

In this subsection we extend $\End^\infty({\mathcal E})$ to an
rapid-decay algebra
$\End^\omega({\mathcal E})$.
Given
a $G$-invariant Dirac-type operator acting on sections of $E$,
we consider the Bismut superconnections $\{A_s\}_{s>0}$ on ${\mathcal E}$.
We compute the $s \rightarrow 0$ limit of the pairing
between the Chern character of $A_s$ and a closed graded trace on
$\Omega^*{\mathcal B}$ that is concentrated on the units $M$.

We now choose
a $G$-invariant vertical Riemannian metric $g^{TZ}$ on the
submersion $\pi \: : \: P \rightarrow M$ and a
$G$-invariant horizontal
distribution $T^H P$.
Given $m \in M$, let $d_m$ denote the corresponding
metric on $Z_m$.  We note that $\{Z_m\}_{m \in M}$ has uniformly
bounded geometry.

Let
$\End^\omega_{{\mathcal B}}
\left( {\mathcal E} \right)$ be the algebra formed by
$G$-invariant operators $K$ as in (\ref{5.20}) whose
integral kernels $K(p^\prime|p) \in
C^\infty( P \times_{M} P; (q^\prime)^* {\mathcal D}_Z \otimes
\Hom((q^\prime)^* E, q^* E))$ are
such that for all $q \in \Z^+$,
\begin{equation}
\sup_{(p^\prime, p) \in P \times_{M} P}
e^{q \: d(p^\prime, p)} \: |K(p^\prime|p)| \: < \: \infty,
\end{equation}
along with the analogous property for the covariant derivatives of $K$.

\begin{proposition}
Let $\rho$ be a linear functional on
$C^\infty_c(M; {\mathcal D})$. Suppose that
the linear functional
$\eta$ on ${\mathcal B}$, defined by
\begin{equation}
\eta(b) \: = \: \rho( b \big|_{M}),
\end{equation}
is a trace on ${\mathcal B}$.
Then $\rho \circ \tau$ is a supertrace on
$\End^\omega_{\mathcal B} ({\mathcal E})$.
\end{proposition}
\begin{pf}
The proof is formally the same as that of
Proposition \ref{prop3}. We omit the details
\end{pf}

Let $\Hom^\omega_{\mathcal B} ({\mathcal E}, \Omega^{*} {\mathcal E})$
be the algebra formed by $G$-invariant operators
$K$ as in (\ref{5.35}) whose integral
kernels
\begin{equation}
K(p^\prime | p) \in C^\infty_c( P \times_M P ;
\Lambda^*((\pi \circ q)^* \tau^*) \otimes (q^\prime)^* {\mathcal D}_Z \otimes
\Hom((q^\prime)^* E, q^* E))
\end{equation}
are
such that for all $q \in \Z^+$,
\begin{equation}
\sup_{(p^\prime, p) \in P \times_{M} P}
e^{q \: d(p^\prime, p)} \: |K(p^\prime|p)| \: < \: \infty,
\end{equation}
along with the analogous property for the covariant derivatives of $K$.

\begin{proposition}
Let $\rho$ be a linear functional on
$C^\infty_c(M; \Lambda^n \tau^* \otimes {\mathcal D})$
Suppose that the linear functional
$\eta$ on $\Omega^n{\mathcal B}$, defined by
\begin{equation}
\eta(\phi) \: = \: \rho( \phi \big|_M),
\end{equation}
is a graded trace on $\Omega^*{\mathcal B}$.
Then $\rho \circ \tau$
is a graded trace on
$\Hom^\omega_{\mathcal B} ({\mathcal E},
\Omega^{*} {\mathcal E})$.
\end{proposition}
\begin{pf}
The proof is formally the same as that of
Proposition \ref{prop4}. We omit the details.
\end{pf}

\begin{proposition}\label{prop6}
Let $\rho$ be a linear functional on
$C^\infty_c(M; \Lambda^n \tau^* \otimes {\mathcal D})$
Suppose that the linear functional
$\eta$ on $\Omega^n{\mathcal B}$, defined by
\begin{equation}
\eta(\phi) \: = \: \rho( \phi \big|_M),
\end{equation}
is a closed graded trace on $\Omega^*{\mathcal B}$.
Then $\rho \circ \tau$
annihilates $[\nabla, K]$ for all $K \in
\Hom^\omega_{\mathcal B} ({\mathcal E},
\Omega^{n-1} {\mathcal E})$.
\end{proposition}
\begin{pf}
The proof is formally the same as that of
Proposition \ref{prop5}. We omit the details.
\end{pf}

Suppose that $Z$ is even-dimensional. Let $E$ be a $G$-invariant
Clifford bundle on $P$ which is equipped with a $G$-invariant
connection. For simplicity of notation, we assume that ${E} \: =
\: S^Z \: \widehat{\otimes} \: \widetilde{V}$, where $S^Z$ is a
vertical spinor bundle and $\widetilde{V}$ is an auxiliary vector
bundle on $P$. More precisely, suppose that the vertical tangent
bundle $TZ$ has a $G$-invariant spin structure.  Let $S^Z$ be the
vertical spinor bundle, a $G$-invariant $\Z_2$-graded Hermitian
vector bundle on $P$. Let $\widetilde{V}$ be another $G$-invariant
$\Z_2$-graded Hermitian vector bundle on $P$ which is equipped
with a $G$-invariant Hermitian connection.
That is, $\widetilde{V}$ is the pullback of a Hermitian
vector bundle $G$ on $P/G$ with a Hermitian connection $\nabla^V$.
Then we put ${E} \: =
\: S^Z \: \widehat{\otimes} \: \widetilde{V}$. The case of general
$G$-invariant Clifford bundles ${E}$ can be treated in a way
completely analogous to what follows.

Let $\nabla^{TZ}$ be the Bismut connection on $TZ$, as constructed
using the horizontal distribution $(d\mu)^{-1}(T^H(P/G))$ on $P$;
see, for example, Berline-Getzler-Vergne
\cite[Proposition 10.2]{Berline-Getzler-Vergne (1992)}.
The $G$-invariance of $\nabla^{TZ}$ and $\nabla^{\widetilde{V}}$
implies that
$\widehat{A}(\nabla^{TZ}) \:
\ch(\nabla^{\widetilde{V}})$ lies in
$C^\infty(P; \Lambda^* (TZ)^* \: \otimes \: \Lambda^*(\pi^* \tau^*))$.

Let $Q \in
\End_{\mathcal B}
\left( {\mathcal E} \right)$
 denote the vertical Dirac-type operator.
From finite-propagation-speed estimates
as in Lott \cite[Proof. of Prop 8]{Lott (1992)},
along with
the bounded geometry
of $\{Z_m\}_{m \in M}$,
for any $s \: > \: 0$ we have
\begin{equation}
e^{- \: s^2 \: Q^2} \in
\End^\omega_{\mathcal B} ({\mathcal E}).
\end{equation}

Let
$A_s \: : \: {\mathcal E} \rightarrow
\Omega^{*} {\mathcal E}$
be the superconnection
\begin{equation}
A_s \: = \: s \: Q \: + \: \nabla^{\mathcal E} \: - \: \frac{1}{4s} \: c(T^P).
\end{equation}
Here
$c(T^P)$ is Clifford multiplication by the curvature $2$-form $T^P$ of
$(d\mu)^{-1}(T^H (P/G))$, restricted to the horizontal vectors
$T^H P$.
We note that the analogous connection term of the Bismut superconnection
\cite[Proposition 10.15]{Berline-Getzler-Vergne (1992)} has an
additional term to make it Hermitian, but in our setting
this term is incorporated into the horizontal differentiation of the
vertical density.
One can use finite-propagation-speed estimates, along with
the bounded geometry
of $\{Z_m\}_{m \in M}$ and the Duhamel expansion as in
\cite[Theorem 9.48]{Berline-Getzler-Vergne (1992)},
to show that we obtain a well-defined element
$e^{- \: (A^s)^2 \: - \: {\mathcal L}_{\widetilde{T}}}$ of
$\Hom^\omega_{\mathcal B} ({\mathcal E},
\Omega^{*} {\mathcal E})$; see
\cite[Theorem 3.1]{Heitsch (1995)} for an analogous statement
when $P \: = \: G \: = \: G_{hol}$.

Let ${\mathcal R}$ be the rescaling operator which, for $p$ even,
multiplies a $p$-form by
$(2 \pi i)^{- \: \frac{p}{2}}$.
Put
\begin{equation}
\ch(A_s) \: = \: {\mathcal R} \left(\tau e^{- \: A_s^2 \: - \:
{\mathcal L}_{\widetilde{T}} }  \right)
\in C^\infty_c(M; \Lambda^* \tau^* \: \otimes \: {\mathcal D}).
\end{equation}

\begin{theorem} \label{theorem1}
Given a linear functional $\rho$
which satisfies the hypotheses of Proposition \ref{prop6},
\begin{equation}
\lim_{s \rightarrow 0}
\rho ( \ch(A_s))
 \: = \:
\rho \left( \int_Z \Phi \: \widehat{A}(\nabla^{TZ}) \:
\ch(\nabla^{\widetilde{V}})  \right).
\end{equation}
\end{theorem}
\begin{pf}
Using Lemmas \ref{lll} and \ref{llll},
$A_s^2 \: + \: {\mathcal L}_{\widetilde{T}}$ is $G$-invariant.
Let $A^\prime_s$ be the corresponding Bismut superconnection
on the foliated manifold $P/G$, a locally-defined differential operator
constructed using
the horizontal distribution $T^H (P/G)$.
By construction, $A_s^2 \: + \: {\mathcal L}_{\widetilde{T}}$
is the pullback
under $\mu$ of $(A^\prime_s)^2$, where we use the identification
$\Lambda^* (\pi^* \tau^*) \: = \: \mu^* \Lambda^* (NF)^*$.
From \cite[Theorem 10.23]{Berline-Getzler-Vergne (1992)}, the
$s \rightarrow 0$ limit of the supertrace of the kernel of
$e^{- \: (A_s^\prime)^2}$, when restricted to the diagonal of
$(P/G) \times (P/G)$, is
$\widehat{A}(\nabla^{TF}) \: \ch(\nabla^{V})$. Then
the
$s \rightarrow 0$ limit of the supertrace of the kernel of
$e^{- \: A_s^2 \: - \:
{\mathcal L}_{\widetilde{T}}}$, when restricted to the diagonal of
$P \times P$, is the pullback under $\mu$ of
$\widehat{A}(\nabla^{TF}) \: \ch(\nabla^{V})$, i.e.
$\widehat{A}(\nabla^{TZ}) \:
\ch(\nabla^{\widetilde{V}})$.
The theorem follows.
\end{pf}
\noindent
{\bf Remark : } If $P \: = \: G \: = \: G_{hol}$ then an
analogue of Theorem \ref{theorem1} appears in
\cite[Theorem 2.1]{Heitsch (1995)}.

If we put
\begin{equation}
G^\prime \: = \: \{ (p_1, p_2) \in P \times P \: : \:
\pi(p_1) \: = \: \pi(p_2)\}/G.
\end{equation}
then $G^\prime$ has the structure of a foliation groupoid, with units
$G^{\prime (0)} \: = \: P/G$. In this way we could reduce from the
case of $G$ acting on $P$ to the case of the foliation groupoid $G^\prime$
acting on itself.  However, doing so would not really simplify any of the
constructions.

\section{Index Theorem} \label{Index Theorem}

In this section we prove the main result of the paper,
Theorem \ref{theorem4}.

\subsection{The index class}

In this subsection we construct the index class
$\Ind(D) \in \KK_0({\mathfrak A})$. We describe its pairing
with a closed graded trace on ${\mathcal B}$. We prove that
the pairing of $\Ind(D)$ with the closed graded trace equals
the pairing of $\ch(A_s)$ with the closed graded trace.

Consider the algebra ${\mathfrak A} \: = \:
\End^\infty_{{\mathcal B}}
\left( {\mathcal E} \right)$.
Let $D \: : \: {\mathcal E}^+ \rightarrow {\mathcal E}^-$ be the
restriction of $Q$ to the positive subspace ${\mathcal E}^+$ of
${\mathcal E}$.
We construct an index projection following
Connes-Moscovici
\cite{Connes-Moscovici (1990)} and Moscovici-Wu
\cite{Moscovici-Wu (1994)}.
Let $u \in C^\infty(\R)$ be an even function such that
$w(x) \: = \: 1 \: - \: x^2 \: u(x)$ is a Schwartz function and
the Fourier transforms of $u$ and $w$ have compact support
\cite[Lemma 2.1]{Moscovici-Wu (1994)}.
Define $\overline{u} \in C^\infty([0, \infty))$ by $\overline{u}(x)
\: = \: u(x^2)$.
Put ${\mathcal P} \: = \: \overline{u}(D^* D) D^*$, which we will think of as a
parametrix for $D$, and put
$S_+ \: = \: I \: - \: {\mathcal P}D$, $S_- \: = \: I \: - \: D{\mathcal P}$.
Consider the operator
\begin{equation}
L \: = \:
\begin{pmatrix}
S_+ & - \: (I + S_+) {\mathcal P} \\
D & S_-
\end{pmatrix},
\end{equation}
with inverse
\begin{equation}
L^{-1} \: = \:
\begin{pmatrix}
S_+ & \: {\mathcal P}(I + S_-) \\
- \: D & S_-
\end{pmatrix}.
\end{equation}
The index projection is defined by
\begin{equation}
p \: = \: L \:
\begin{pmatrix}
I & 0 \\
0 & 0
\end{pmatrix} \:
L^{-1} \: = \:
\begin{pmatrix}
S_+^2 & S_+ (I + S_+) {\mathcal P} \\
S_- D & I - S_-^2
\end{pmatrix}.
\end{equation}
Put
\begin{equation}
p_0 \: = \:
\begin{pmatrix}
0 & 0 \\
0 & I
\end{pmatrix}.
\end{equation}
By definition, the index of $D$ is
\begin{equation}
\Ind(D) \: = \: [p \: - \: p_0] \in \KK_0({\mathfrak A}).
\end{equation}

As $Q$ is $G$-invariant, the operator $l$ of (\ref{5.5})
commutes with $p$, and (\ref{5.6}) holds for
$\xi \in \Image(p)$.
If $\rho$ is a linear functional
which satisfies the hypotheses of Proposition \ref{prop5}, define
the pairing of $\rho$ with $\Ind (D)$ by
\begin{align} \label{ch1}
\langle \ch (\Ind (D))   , \rho \rangle \: = \:
& (2 \pi i )^{- \: deg(\rho)/2} \\
& \rho \left(
 \tau \left( p \: e^{- (p \circ \nabla^{\mathcal E} \circ p)^2 \: - \:
{\mathcal L}_{\widetilde{T}}} \: p  \: - \:
p_0 \: e^{- (p_0 \circ \nabla^{\mathcal E} \circ p_0)^2
\: - \:
{\mathcal L}_{\widetilde{T}}} \: p_0 \right) \right), \notag
\end{align}
where we have extended the ungraded trace $\tau$ in the obvious
way to act on $(2 \times 2)$-matrices.
(See \cite[Section 5]{Gorokhovsky-Lott (2003)} for the justification
of the definition.)

\begin{theorem} \label{theorem2}
For all $s > 0$,
\begin{equation}\label{eqn1}
\langle \ch(\Ind (D)), \rho \rangle \: = \:
\rho ( \ch(A_s)).
\end{equation}
\end{theorem}
\begin{pf}
The proof follows the lines of the proof of
\cite[Proposition 4 and Theorem 3]{Gorokhovsky-Lott (2003)}, to
which we refer for details. We only
present the main idea.
Put
\begin{equation}
\nabla^\prime \: = \: \left(
\begin{pmatrix}
I & 0 \\
0 & 0
\end{pmatrix}
 \: L^{-1} \circ \nabla^{\mathcal E} \circ L \:
\begin{pmatrix}
I & 0 \\
0 & 0
\end{pmatrix}
\right)\: + \: \left(
\begin{pmatrix}
0 & 0 \\
0 & I
\end{pmatrix}
 \: \nabla^{\mathcal E} \:
\begin{pmatrix}
0 & 0 \\
0 & I
\end{pmatrix} \right).
\end{equation}
Then one can show algebraically that
\begin{equation}
\langle \ch(\Ind (D)), \rho \rangle \: = \:
\rho \left( {\mathcal R} \: \tau e^{- \: (\nabla^\prime)^2 \: - \:
{\mathcal L}_{\widetilde{T}}} \right),
\end{equation}
where the $\tau$ on the right-hand-side is now a graded trace.
Next, one shows that
\begin{equation}
\rho \left( {\mathcal R} \: \tau e^{- \: (\nabla^\prime)^2 \: - \:
{\mathcal L}_{\widetilde{T}}} \right)
\: = \:
\rho ( \ch(A_s))
\end{equation}
by performing a homotopy
from $\nabla^\prime$ to $A_s$, from which the
theorem follows.
\end{pf}

\subsection{Construction of $\omega_\rho$}\label{subsection6.4}

In this subsection we
construct the universal class $\omega_\rho \in
\HH^*(BG; o)$.
We express $\rho(\ch(A_s))$ as an integral
involving the pullback of $\omega_\rho$.

Put $V \: = \: \widetilde{V}/G$, a
Hermitian vector bundle on $P/G$ with a compatible connection $\nabla^V$.

Let
$o(\tau)$ be the orientation bundle of $\tau$, a flat real line bundle on
$M$.
Let $\rho$ satisfy the hypotheses of Proposition \ref{prop5}.
By duality, $\rho$
corresponds to a closed distributional form
$* \rho \in \Omega^{\dim(M)-n}(M; o(\tau))$.

Let $EG$ denote the bar construction of a universal space on which
$G$ acts freely. That is, put
\begin{equation}
G^{(n)} \: = \: \{(g_1, \ldots, g_n) \: : \: s(g_1) \: = \: r(g_2),
\ldots, s(g_{n-1}) \: = \: r(g_n) \}.
\end{equation}
Then
$EG$ is the geometric realization of a
simplicial manifold given by $E_nG \: = \: G^{(n+1)}$, with
face maps
\begin{equation} \label{5.47}
d_i(g_0, \ldots, g_n) \: = \:
\begin{cases}
(g_1, \ldots, g_n) &\text{ if } i = 0,\\
(g_0, \ldots, g_{i-1} g_i, \ldots, g_n) &
\text{ if } 1 \: \le \:  i \: \le \: n
\end{cases}
\end{equation}
and degeneracy maps
\begin{equation} \label{5.48}
s_i(g_0, \ldots, g_n) \: = \:
(g_0, \ldots, g_{i}, 1, g_{i+1}, \ldots g_n),
\: \: \: \: \: 0 \: \le \: i \: \le \: n.
\end{equation}
Here $1$ denotes a unit.
The action of $G$ on $EG$ is induced from the action on $E_nG$ given
by $(g_0, \ldots, g_n) \: g \: = \:
(g_0, \ldots, g_n g)$. Let $BG$ be the quotient space.
Define $\pi^\prime \: : \: EG \rightarrow M$ as the extension of
$(g_0, \ldots, g_n) \rightarrow  s(g_n)$.
Put
$\Omega^{n_1, n_2}(EG) \: = \: \Omega^{n_1}(G^{(n_2+1)})$ and
$\Omega^{n_1, n_2}(BG) \: = \: \left( \Omega^{n_1, n_2}(EG) \right)^G$.
Let $\Omega^*(BG)$ be the total complex of $\Omega^{*,*}(BG)$.
Here the forms on $G^{(n_2+1)}$
can be either smooth or distributional, depending on
the context. We will speak correspondingly of smooth or distributional
elements of $\Omega^*(BG)$. In either case,
the cohomology of $\Omega^*(BG)$ equals
$\HH^*(BG; \R)$. There is a similar
discussion for twistings by a local system.

The action of $G$ on $P$
is classified by a continuous $G$-equivariant map
$\widehat{\nu} \: : \: P \rightarrow EG$.
Let $\nu \: : \: P/G \rightarrow BG$ be the $G$-quotient of
$\widehat{\nu}$.
There are commutative
diagrams
\begin{equation}
\begin{CD}
P @>\widehat{\nu}>> EG \\
@V{\pi}VV @V{\pi^\prime}VV \\
M
@>\Id.>>
M
\end{CD}
\end{equation}
and
\begin{equation}
\begin{CD}
P @>\widehat{\nu}>> EG \\
@V VV @V VV \\
P/G
@>\nu>>
BG.
\end{CD}
\end{equation}
As $P/G$ is compact, we may assume that $\nu$ is
Lipschitz.

Consider $(\pi^\prime)^* (*\rho) \in
\Omega^{*}(EG; (\pi^\prime)^* o(\tau))$, a closed
distributional
form on $EG$. Let $o$ be the $G$-quotient of $(\pi^\prime)^* o(\tau)$,
a flat real line bundle on $BG$. Then
$(\pi^\prime)^* (*\rho)$ pulls back from a closed distributional
form in $\Omega^*(BG; o)$, which represents a class in
$\HH^*(BG; o)$.
Let $\omega_\rho \in \Omega^*(BG; o)$ be a closed
smooth
form representing the same cohomology class. Let $\widehat{\omega}_\rho \in
\Omega^{*}(EG; (\pi^\prime)^* o(\tau))$
be its pullback to $EG$. As $\nu$ is Lipschitz,
$\nu^* \omega_\rho$ is an $L^\infty$-form on $P/G$.

\begin{theorem} \label{theorem3}
\begin{equation}
\rho \left( \int_Z \Phi \: \widehat{A}(\nabla^{TZ}) \:
\ch(\nabla^{\widetilde{V}})  \right) \: = \:
\int_{P/G} \widehat{A}({TF}) \: \ch(V) \: \nu^* \omega_\rho.
\end{equation}
\end{theorem}
\begin{pf}
Let $* \left( \Phi \: \widehat{A}(\nabla^{TZ}) \:
\ch(\nabla^{\widetilde{V}}) \right)$ be the dual of
$\Phi \: \widehat{A}(\nabla^{TZ}) \:
\ch(\nabla^{\widetilde{V}})$. We will think of
$* \left( \Phi \: \widehat{A}(\nabla^{TZ}) \:
\ch(\nabla^{\widetilde{V}}) \right)$ as a cycle on $P$
and $(\pi^\prime)^* (*\rho)$ as a cocycle on $EG$. Then
\begin{align}
\rho \left( \int_Z \Phi \: \widehat{A}(\nabla^{TZ}) \:
\ch(\nabla^{\widetilde{V}})  \right) \:
& = \: \langle \pi_* \left(
* \left( \Phi \: \widehat{A}(\nabla^{TZ}) \:
\ch(\nabla^{\widetilde{V}}) \right)
\right),  *\rho \rangle_{M} \\
& = \: \langle
* \left( \Phi \: \widehat{A}(\nabla^{TZ}) \:
\ch(\nabla^{\widetilde{V}}) \right),
\pi^* (*\rho) \rangle_{P} \notag \\
& = \: \langle \widehat{\nu}_* \left(
* \left( \Phi \: \widehat{A}(\nabla^{TZ}) \:
\ch(\nabla^{\widetilde{V}}) \right)
\right), (\pi^\prime)^* (*\rho) \rangle_{EG} \notag \\
& = \: \langle \widehat{\nu}_* \left(
* \left( \Phi \: \widehat{A}(\nabla^{TZ}) \:
\ch(\nabla^{\widetilde{V}}) \right)
\right), \widehat{\omega}_\rho \rangle_{EG} \notag \\
& = \: \langle
* \left( \Phi \: \widehat{A}(\nabla^{TZ}) \:
\ch(\nabla^{\widetilde{V}}) \right),
\widehat{\nu}^*
\widehat{\omega}_\rho \rangle_{P} \notag \\
& = \:
\int_P \Phi \: \widehat{A}(\nabla^{TZ}) \:
\ch(\nabla^{\widetilde{V}}) \: \widehat{\nu}^* \widehat{\omega}_\rho
\notag \\
& = \:
\int_{P/G} \widehat{A}(TF) \:
\ch(V) \: {\nu}^* {\omega}_\rho. \notag
\end{align}
\end{pf}
\noindent
{\bf Remark : } If one were willing to work with orbifolds $P/G$
instead of manifolds then one could extend Theorem \ref{theorem3} to
general proper cocompact actions, with
$\omega_\rho \in \HH^*(\underline{B}G; o)$ being a
cohomology class
on the classifying space for proper $G$-actions.

\subsection{Proof of index theorem}

\begin{theorem} \label{theorem4}
If $G$ acts freely, properly discontinuously and cocompactly on
$P$  and $\rho$ satisfies the hypotheses of Proposition \ref{prop6} then
\begin{equation}
\langle \ch(\Ind D), \rho \rangle\: = \: \int_{P/G}
\widehat{A}({TF}) \: \ch(V) \: \nu^* \omega_\rho.
\end{equation}
\end{theorem}
\begin{pf}
If $Z$ is even-dimensional then the claim follows from Theorems
\ref{theorem1}, \ref{theorem2} and \ref{theorem3}.
If $Z$ is odd-dimensional
then one can reduce to the even-dimensional case by a standard trick
involving taking the product with a circle.
\end{pf}
{\bf Example 8 :} Suppose that $(M, {\mathcal F})$ is a closed
foliated manifold.
Take $P \: = \: G \: = \: G_{hol}$. Let $\mu$ be a holonomy-invariant
transverse measure for ${\mathcal F}$. Take $\rho$ as in Example 6.  Then
Theorem \ref{theorem4} reduces to Connes' $L^2$-foliation index theorem
\cite[Section I.5.$\gamma$, Theorem 7]{Connes (1994)}
\begin{equation}
\langle \Ind D, \rho \rangle\: = \: \langle
\widehat{A}({TF}) \: \ch(V), RS_\mu \rangle,
\end{equation}
where $RS_\mu$ is the Ruelle-Sullivan current associated to $\mu$
\cite[Section I.5.$\beta$]{Connes (1994)}. \\ \\
{\bf Example 9 :} Let $(M, {\mathcal F})$ be a closed manifold equipped
with a codimension-$q$ foliation.
Take $P \: = \: G \: = \: G_{hol}$. 
Let $\HH^*(\Tr {\mathcal F})$ denote the Haefliger cohomology of
$(M, {\mathcal F})$ \cite{Haefliger (1980)}. Recall that there is a linear map
$\int_{\mathcal F} \: : \: \HH^*(M) \rightarrow \HH^{*-n+q}(\Tr
{\mathcal F})$. 
Let $c$ be a closed holonomy-invariant
transverse current for ${\mathcal F}$. Take $\rho$ as in Example 7.
Then
Theorem \ref{theorem4} becomes
\begin{equation}
\langle \ch(\Ind D), \rho \rangle\: = \: \langle
\int_{\mathcal F} \widehat{A}({TF}) \: \ch(V), c \rangle.
\end{equation}
This is a consequence of the Connes-Skandalis
foliation index theorem, along with the
result of Connes that $\rho$ gives a higher trace on the
reduced foliation $C^*$-algebra; see
\cite{Benameur-Heitsch (2005),Connes (1986),Connes-Skandalis (1984)}. \\ \\
{\bf Example 10 : } Let $M$ be a closed oriented $n$-dimensional
manifold. Let $G = M$ be the
groupoid that just consists of units.
Let $P$ be a closed manifold that is the total space
of an oriented fiber bundle $\pi \: : \: P \rightarrow M$ with fiber $Z$.
Let $c$ be a closed current on $M$
with homology class $[c] \in \HH_*(M; \C)$.
With $* \: : \: 
\HH_*(M; \C) \rightarrow \HH^{n-*}(M; \C)$ being the Poincar\'e
isomorphism,
Theorem \ref{theorem4} becomes
\begin{equation}
\langle \ch(\Ind D), c \rangle\: = \: \int_P
\widehat{A}({TZ}) \: \ch(V) \: \pi^*(*[c])).
\end{equation}
This is a consequence of the Atiyah-Singer families index theorem
\cite{Atiyah-Singer (1971)}, as the right-hand-side equals
$\langle \int_Z
\widehat{A}({TZ}) \: \ch(V), c \rangle$.

{\bf Example 11 : } Let $G$ be a discrete group that acts
freely, properly discontinuously and cocompactly on a manifold $P$.
As its space of units $M$ is a point, let $\rho$ be the identity map
$C^\infty(M) \rightarrow \C$. Then
Theorem \ref{theorem4} reduces to Atiyah's $L^2$-index theorem
\cite{Atiyah (1976)}
\begin{equation}
\langle \Ind D, \rho \rangle\: = \:
\int_{P/G} \widehat{A}({TP/G}) \: \ch(V).
\end{equation}

\appendix
\section{Appendix}

This is an addendum to \cite{Gorokhovsky-Lott (2003)}, in
which we use finite propagation speed methods to improve
\cite[Theorem 3]{Gorokhovsky-Lott (2003)}. In
the improved version we allow $\eta$ to be a closed
graded trace on $\Omega^*(B, \C \Gamma)$, as opposed to
$\Omega^*(B, {\mathcal B}^\omega)$. There is a similar
improvement of \cite[Theorem 6]{Gorokhovsky-Lott (2003)}.

We will follow the notation of \cite{Gorokhovsky-Lott (2003)}.

\subsection{Finite propagation speed}\label{subsection6.1}

Let $f \in C^\infty_c(\R)$ be a smooth even function with support
in $[-\epsilon, \epsilon]$. Put
\begin{equation}
\widehat{f}(y) \: = \: \int_{\R} f(x) \: \cos(xy) \: dx,
\end{equation}
a smooth even function.
With $A_s$ as in \cite[(4.7)]{Gorokhovsky-Lott (2003)},
put
\begin{equation} \label{fourier}
\widehat{f} ( A_s ) \: = \:
\int_{\R} f(x) \: \cos \left( x \: A_s \right) \: dx.
\end{equation}

Let us describe $\cos \left( x \:
A_s \right)$ explicitly, using the fact
that it satisfies
\begin{equation} \label{wave}
\left( \partial_x^2 \: + \: A_s^2
\right) \:  \cos \left( x \: A_s
 \right) \: = \: 0.
\end{equation}
Write
$A_s^2
\: = \: s^2 Q^2 \: + \: X$. We first consider a solution $u(\cdot, x)$
of the inhomogeneous wave equation
\begin{equation}
\left( \partial_x^2 \: + \: s^2 Q^2 \right) \: u \: = \: f
\end{equation}
with initial conditions $u(\cdot, 0) = u_0(\cdot)$ and
$u_x(\cdot, 0) = 0$. Then $u(\cdot, x)$ is given by
\begin{equation}
u(x) \: = \: \cos(xsQ) u_0 \: + \: \int_0^x
\frac{\sin((x-v)sQ)}{sQ} \: f(v) \: dv.
\end{equation}
Putting $f \: = \: - \: X u$ and iterating, we obtain an expansion
of $\cos \left( x \: A_s \right)$ of the form
\begin{equation}
\cos \left( x \: A_s \right) \: = \:
\cos(xsQ) \: - \: \int_0^x
\frac{\sin((x-v)sQ)}{sQ} \: X \: \cos(vsQ) \: dv \: + \: \ldots
\end{equation}
Because $X$ has positive form degree, there is no problem with
the convergence of the series.

From finite propagation speed results, we know that
$\cos(xsQ)$
has a Schwartz kernel $\cos(xsQ)(p^\prime | p)$ with support
on $\{(p^\prime, p) \: : \: d(p^\prime, p) \: \le \: xs\}$, and
similarly for $\frac{\sin(xsQ)}{sQ}$; see Taylor
\cite[Chapter 4.4]{Taylor (1981)}. Using the compactness of
$h$, it follows that the $(m,n)$-component
$\widehat{f} ( A_s )_{(m,n)}$ lies in
$\Hom^\infty_{C^\infty_c(B) \rtimes \Gamma}
(C^\infty_c(\widehat{M}; \widehat{E}), \Omega^{m,n}(B, \C \Gamma)
\otimes_{C^\infty_c(B) \rtimes \Gamma}
C^\infty_c(\widehat{M}; \widehat{E}))$.

Finally, define
$\ch_{\widehat{f}}(A_s) \in \Omega^*(B, \C \Gamma)_{ab}$ by
\begin{equation}
\ch_{\widehat{f}}(A_s) \: = \: {\mathcal R} \Tr_{s, \langle e \rangle}
\widehat{f}(A_s).
\end{equation}

\subsection{Index Pairing}\label{subsection6.2}

In this subsection we show that for all $s > 0$ and all
closed graded traces $\eta$ on $\Omega^*(B, \C \Gamma)$,
$\langle \ch_{\widehat{f}}(A_s), {\eta} \rangle
\: = \: \langle \widehat{f}(\Ind (D)), {\eta} \rangle$.
The method of proof is essentially
the same as that of \cite[Section 5]{Gorokhovsky-Lott (2003)},
which in turn was inspired by Nistor
\cite{Nistor (1997)}.

In analogy to \cite[Section 5.3]{Gorokhovsky-Lott (2003)}, put
${\mathcal E} \: = \: C^\infty_c(\widehat{M}; \widehat{E})$ and
$\widetilde{\mathfrak A} \: = \:
\End^\infty_{C^\infty_c(B) \rtimes \Gamma}
\left(  C^\infty_c(\widehat{M}; \widehat{E}) \right)$.
Let $D \: : \: {\mathcal E}^+ \rightarrow {\mathcal E}^-$ be the
restriction of $Q$ to the positive subspace ${\mathcal E}^+$ of
${\mathcal E}$.
We construct an index projection following
\cite{Connes-Moscovici (1990)} and
\cite{Moscovici-Wu (1994)}.
Let $u \in C^\infty(\R)$ be an even function such that
$w(x) \: = \: 1 \: - \: x^2 \: u(x)$ is a Schwartz function and
the Fourier transforms of $u$ and $w$ have compact support
\cite[Lemma 2.1]{Moscovici-Wu (1994)}.
Define $\overline{u} \in C^\infty([0, \infty))$ by $\overline{u}(x)
\: = \: u(x^2)$.
Put ${\mathcal P} \: = \: \overline{u}(D^* D) D^*$, which we will think of as a
parametrix for $D$, and put
$S_+ \: = \: I \: - \: {\mathcal P}D$, $S_- \: = \: I \: - \: D{\mathcal P}$.
Consider the operator
\begin{equation}
L \: = \:
\begin{pmatrix}
S_+ & - \: (I + S_+) {\mathcal P} \\
D & S_-
\end{pmatrix},
\end{equation}
with inverse
\begin{equation}
L^{-1} \: = \:
\begin{pmatrix}
S_+ & \: {\mathcal P}(I + S_-) \\
- \: D & S_-
\end{pmatrix}.
\end{equation}
The index projection is defined by
\begin{equation}
p \: = \: L \:
\begin{pmatrix}
I & 0 \\
0 & 0
\end{pmatrix} \:
L^{-1} \: = \:
\begin{pmatrix}
S_+^2 & S_+ (I + S_+) {\mathcal P} \\
S_- D & I - S_-^2
\end{pmatrix}.
\end{equation}
Put
\begin{equation}
p_0 \: = \:
\begin{pmatrix}
0 & 0 \\
0 & I
\end{pmatrix}.
\end{equation}
By definition, the index of $D$ is
\begin{equation}
\Ind(D) \: = \: [p \: - \: p_0] \in \KK_0(\widetilde{\mathfrak A}).
\end{equation}

Put $\widetilde{\Omega}^* \: = \:
\Hom^\infty_{C^\infty_c(B) \rtimes \Gamma}
(C^\infty_c(\widehat{M}; \widehat{E}), \Omega^{*}(B, \C \Gamma)
\otimes_{C^\infty_c(B) \rtimes \Gamma}
C^\infty_c(\widehat{M}; \widehat{E}))$, a graded algebra
with derivation
$\nabla \: = \: \nabla^{(1,0)} \: + \: \nabla^{(0,1)}$.
If $\eta$ is a closed graded trace on $\Omega^*(B, \C \Gamma)$,
define
the pairing of $\eta$ with $\Ind (D)$ by
\begin{equation} \label{ch2}
\langle \widehat{f} (\Ind (D))   , {\eta} \rangle \: = \:
(2 \pi i )^{- \: deg(\eta)/2} \:
\langle \Tr_{\langle e \rangle}
\left(\widehat{f}(p \circ \nabla \circ p) \: - \:
\widehat{f}( p_0 \circ \nabla \circ p_0) \right), \eta \rangle.
\end{equation}
(See \cite[Section 5]{Gorokhovsky-Lott (2003)} for the justification
of the definition.)

\begin{theorem} \label{tthm}
For all $s > 0$,
\begin{equation}\label{eqn2}
\langle \ch_{\widehat{f}}(A_s), {\eta} \rangle
\: = \: \langle \widehat{f}(\Ind (D)), {\eta} \rangle.
\end{equation}
\end{theorem}
\begin{pf}
The proof follows the lines of the proof of
\cite[Proposition 4 and Theorem 3]{Gorokhovsky-Lott (2003)}, to
which we refer for details. We only
present the main idea.
Put
\begin{equation}
\nabla^\prime \: = \: \left(
\begin{pmatrix}
I & 0 \\
0 & 0
\end{pmatrix}
 \: L^{-1} \circ \nabla \circ L \:
\begin{pmatrix}
I & 0 \\
0 & 0
\end{pmatrix}
\right)\: + \: \left(
\begin{pmatrix}
0 & 0 \\
0 & I
\end{pmatrix}
 \: \nabla \:
\begin{pmatrix}
0 & 0 \\
0 & I
\end{pmatrix} \right).
\end{equation}
Then one can show algebraically that
\begin{equation}
\langle \widehat{f}(\Ind (D)), \eta \rangle \: = \:
\langle {\mathcal R} \: \Tr_{s, \langle e \rangle}
\widehat{f} (\nabla^\prime), \eta \rangle.
\end{equation}
Next, one shows that
\begin{equation}
\langle {\mathcal R} \: \Tr_{s, \langle e \rangle}
\widehat{f} (\nabla^\prime), \eta \rangle \: = \:
\langle \ch_{\widehat{f}}(A_s), \eta \rangle
\end{equation}
by performing a homotopy
from $\nabla^\prime$ to $A_s$, from which the
theorem follows. The argument is the same as in
the proof of \cite[Proposition 4]{Gorokhovsky-Lott (2003)}. We
refer to \cite{Gorokhovsky-Lott (2003)}, and will
only indicate the necessary modifications of the equations in
\cite[Section 5.2]{Gorokhovsky-Lott (2003)}.

As in \cite[(5.20)]{Gorokhovsky-Lott (2003)}, for $t \in [0,1]$ put
\begin{equation} \label{4.21}
A(t) \: = \:
\begin{pmatrix}
(\nabla^\prime)^+ & t \: D^* \\
t \: D &  (\nabla^\prime)^-
\end{pmatrix}.
\end{equation}
The analog of \cite[(5.26)]{Gorokhovsky-Lott (2003)} is
\begin{align} \label{add3}
&
\cos \left( x \: A(t) \right) \equiv
\\
&
\begin{pmatrix}
\cos \left( x \sqrt{((\nabla^\prime)^+)^2 \: + \: t^2 \: D^* D}
 \right) &
{\mathcal Z}  \\
0 &
D \:
\cos \left(x \sqrt{ ((\nabla^\prime)^+)^2 \: + \: t^2 \: D^* D}
\right)
\: {\mathcal P}
\end{pmatrix}, \notag
\end{align}
where
\begin{align} \label{add4}
{\mathcal Z} \: = \: - \: & \int_0^x
\frac{\sin \left( (x-v) \sqrt{((\nabla^\prime)^+)^2 \: + \: t^2 \: D^* D}
\right)}{\sqrt{((\nabla^\prime)^+)^2 \: + \: t^2 \: D^* D}} \\
& \left( t \: [(\nabla^\prime)^-, D^*] \: + \: t ((\nabla^\prime)^+ \: - \:
(\nabla^\prime)^-) D^* \right)
\cos \left( v \:
\sqrt{ (\nabla^-)^2 \: + \: t^2 \: D D^*} \right) \: dv \notag
\end{align}
and the left-hand-side of (\ref{add3}) is to be multiplied by
$f$ and then integrated.
As in \cite[(5.30)]{Gorokhovsky-Lott (2003)},
\begin{equation}
\frac{dA}{dt} \: = \: \begin{pmatrix}
0 & D^* \\
D & 0
\end{pmatrix}
\end{equation}
The analog of \cite[(5.31)]{Gorokhovsky-Lott (2003)} is
\begin{align}
& \Tr_s \left( \frac{dA}{dt} \:
\begin{pmatrix}
\cos \left( x \sqrt{((\nabla^\prime)^+)^2 \: + \: t^2 \: D^* D}
\right) &
{\mathcal Z}  \\
0 &
D \:
\cos \left(x \sqrt{ ((\nabla^\prime)^+)^2 \: + \: t^2 \: D^* D
 }\right)
\: {\mathcal P}
\end{pmatrix} \right) \\
& = \:
- \: \Tr \left( D \: {\mathcal Z} \right) \: = \notag \\
& t \: \Tr \left( D \:
\int_0^x
\frac{\sin \left( (x-v) \sqrt{((\nabla^\prime)^+)^2 \: + \: t^2 \: D^* D}
\right)}{\sqrt{((\nabla^\prime)^+)^2 \: + \: t^2 \: D^* D
}}
\: \left( [(\nabla^\prime)^-, D^*] \: + \: ((\nabla^\prime)^+ \: - \:
(\nabla^\prime)^-) D^* \right)
 \right. \notag \\
& \left.
\cos \left( v \:
\sqrt{ (\nabla^-)^2 \: + \: t^2 \: D D^*
} \right) \right) \: dv.
 \notag
\end{align}
The analog of \cite[(5.32)]{Gorokhovsky-Lott (2003)} is
\begin{align}
& D \:
\int_0^x
\frac{\sin \left( (x-v) \sqrt{((\nabla^\prime)^+)^2 \: + \: t^2 \: D^* D
 }
\right)}{\sqrt{((\nabla^\prime)^+)^2 \: + \: t^2 \: D^* D
}}
\: \left( [(\nabla^\prime)^-, D^*] \: + \: ((\nabla^\prime)^+ \: - \:
(\nabla^\prime)^-) D^* \right)
\\
& \cos \left( v \:
\sqrt{ (\nabla^-)^2 \: + \: t^2 \: D D^*
} \right) \: dv \: \equiv \: \notag \\
&
\int_0^x
\frac{\sin \left( (x-v) \sqrt{(\nabla^-)^2 \: + \: t^2 \: D D^*}
\right)}{\sqrt{(\nabla^-)^2 \: + \: t^2 \: DD^*}} \: D
\: \left( [(\nabla^\prime)^-, D^*] \: + \: ((\nabla^\prime)^+ \: - \:
(\nabla^\prime)^-) D^* \right) \notag
\\
& \cos \left( v \:
\sqrt{ (\nabla^-)^2 \: + \: t^2 \: D D^*}
\right)  \: dv \: \equiv \: \notag \\
& \int_0^x
\frac{\sin \left( (x-v) \sqrt{(\nabla^-)^2 \: + \: t^2 \: D D^*}
\right)}{\sqrt{(\nabla^-)^2 \: + \: t^2 \: D D^*
}} \: [\nabla^-, DD^*] \: \notag \\
& \cos \left( v \:
\sqrt{ (\nabla^-)^2 \: + \: t^2 \: D D^*}
\right) \: dv. \notag
\end{align}
The analog of \cite[(5.33)]{Gorokhovsky-Lott (2003)} is
\begin{align}
& \Tr \left(
\int_0^x
\frac{\sin \left( (x-v) \sqrt{(\nabla^-)^2 \: + \: t^2 \: DD^*}
\right)}{\sqrt{(\nabla^-)^2 \: + \: t^2 \: DD^*
}} \: [\nabla^-, DD^*] \right. \\
& \left. \cos \left( v \:
\sqrt{ (\nabla^-)^2 \: + \: t^2 \: D D^*}
\right) \: dv \right) \: =  \notag \\
& - \: t^{-2} \: d \: \Tr \left( \cos \left( x \:
\sqrt{ (\nabla^-)^2 \: + \: t^2 \: D D^*}
\right) \right). \notag
\end{align}
The rest of the proof is as in
\cite[Proof of Proposition 4]{Gorokhovsky-Lott (2003)}.
\end{pf}

We define $\langle \ch(\Ind (D)), \eta \rangle$ by formally taking
$\widehat{f}(z) \: = \: e^{-z^2}$ in (\ref{ch2}). This makes perfect sense,
given that $\eta$ acts on elements of a fixed degree.

\begin{corollary} \label{cor}
a. The left-hand-side of (\ref{eqn2}) only depends on $f$ through the
derivative $\widehat{f}^{(deg(\eta))}(0)$.\\
b. If
$\widehat{f}^{(deg(\eta))}(0) \: = \:
\frac{d^{deg(\eta)} e^{-z^2}}{d^{deg(\eta)} z}  \Big|_{z=0}$ then
\begin{equation}
\langle \ch(\Ind (D)), \eta \rangle \: = \:
\langle \ch_{\widehat{f}}(A_s), \eta \rangle.
\end{equation}
\end{corollary}
\begin{pf}
a. From (\ref{ch2}), the right-hand-side of (\ref{eqn2})
only depends on $f$ through the
derivative $\widehat{f}^{(deg(\eta))}(0)$. From Theorem \ref{tthm}, the
same must be true of the left-hand-side. \\
b. If
$\widehat{f}^{(deg(\eta))}(0) \: = \:
\frac{d^{deg(\eta)} e^{-z^2}}{d^{deg(\eta)} z}  \Big|_{z=0}$
then $\widehat{f}$ has the same relevant term
in its Taylor expansion as the function $z \rightarrow
e^{ - z^2}$, from which the corollary follows.
\end{pf}

\subsection{Pairing of the Chern character of the index with
general closed graded traces} \label{subsection6.3}

In this subsection we prove a formula for the pairing of
the Chern character of the index with a closed graded trace $\eta$
on $\Omega^*(B, \C \Gamma)$.
The idea is to approximate the Gaussian function,
which was previously used in
forming the superconnection Chern character,
by an appropriate function $\widehat{f}$.

\begin{theorem}
Given a closed graded trace $\eta$ on $\Omega^*(B, \C \Gamma)$,
\begin{equation} \label{toshow}
\langle \ch(\Ind (D)), \eta \rangle
 \: = \: \langle
 \int_Z \Phi \: \widehat{A}(\nabla^{TZ}) \:
\ch(\nabla^{\widetilde{V}}) \:
e^{- \: \frac{\nabla_{can}^2}{2 \pi i}},
\eta \rangle.
\end{equation}
\end{theorem}
\begin{pf}
Choose an even function $f \in C^\infty_c(\R)$ so that $\widehat{f}$
satisfies the hypothesis of Corollary \ref{cor}.b.
By Corollary \ref{cor}, it suffices to compute
\begin{equation}
\lim_{s \rightarrow 0}
\langle \ch_{\widehat{f}}(A_s), \eta \rangle.
\end{equation}
With reference to (\ref{fourier}), the local supertrace
$\tr_s \cos \left( x \: A_s \right)(p, p)$ exists as a distribution
in $x$. The singularities near $x = 0$ of the distribution have coefficients
that are the same, up to constants, as the leading terms in
the $x$-expansion of
$\tr_s e^{-x^2 \: A_s^2}(p, p)$;
see, for example, Sandoval
\cite{Sandoval (1999)} for the analogous statement for $\cos(xsQ)$.
As in \cite[Lemma 10.22]{Berline-Getzler-Vergne (1992)}, these are
the terms that enter into the local index computation.
Now $\cos \left( x \: A_s \right)$
satisfies (\ref{wave}),
in analogy to the fact that
$e^{- t \: A_s^2}$ satisfies the heat equation
\begin{equation}
\left( \partial_t \: + \: A_s^2 \right) \: e^{- \: t \: A_s^2} \: = \: 0.
\end{equation}
We can perform a Getzler rescaling as in the proof of
\cite[Theorem 2]{Gorokhovsky-Lott (2003)},
to see that for the purposes of computing the local index, we can
effectively replace the $A_s^2$-term in the differential operator of
(\ref{wave}) by
\cite[(4.12)]{Gorokhovsky-Lott (2003)}. Thus we are reduced to considering
the wave operator of the harmonic oscillator Hamiltonian.  The
rest of the proof of the theorem can
in principle be carried out in a way similar to that
of \cite[Theorem 2]{Gorokhovsky-Lott (2003)}.
However, we can shortcut the calculations
by noting that Corollary \ref{cor}, along with the choice of $f$,
implies that the result of the local calculation must be the same as
$\lim_{s \rightarrow 0} \langle \ch(A_s), \eta \rangle$,
which was already calculated in
\cite[Theorem 2]{Gorokhovsky-Lott (2003)}.
\end{pf}

\end{document}